\documentclass[11pt]{article}
\usepackage{placeins,graphicx,diagbox }

\usepackage{graphicx,bm}
\usepackage{amsfonts,amsmath,fullpage,bbm}
\usepackage{amssymb,multirow,verbatim}
\usepackage{acronym,wrapfig,plain,mathrsfs,enumerate,relsize,color}
\newtheorem{algorithm}{Algorithm}
\usepackage{algorithm}
\newtheorem{assumption}{Assumption}
\usepackage{subfig}
\usepackage{algorithmic}
\usepackage{pifont}

\newcommand{\us}[1]{{\color{black}#1}}

\usepackage{wrapfig} 
\usepackage{times} 
\usepackage{amsmath,mathtools,float}
\usepackage{tikz,pgfplots}
\usepackage{booktabs,caption}
\usepackage[flushleft]{threeparttable}
\usetikzlibrary{arrows,shapes,snakes,shadows,positioning,automata,patterns}
\usetikzlibrary{trees,decorations.pathmorphing,decorations.markings}
\usepackage{amsmath}
\usepackage{multirow}
\usepackage{amssymb,latexsym}  
\usepackage{dsfont} 

\usepackage{amsmath,algorithm} 
\usepackage{amssymb}  
\usepackage{color}
\usepackage{cite}

\newcommand{\blue}[1]{{\color{black}#1}}
\newcommand{\uvs}[1]{{\color{black}#1}}
\newcommand{\red}[1]{{\color{black}#1}}

\newcommand{\epsilonbar}{{\bar \epsilon}}

\def\Real{\mathbb{R}}

\def\argmin{\mathop{\rm argmin}}

\newtheorem{theorem}{Theorem}
\newtheorem{lemma}{Lemma}

\newtheorem{proposition}{Proposition}
\newtheorem{corollary}{Corollary}
\newtheorem{remark}{Remark}

\def\us#1{{{\color{black}#1}}}
\usepackage{verbatim}
\usepackage{ulem}
\begin{document}

\title{Asynchronous  Variance-reduced Block  Schemes    for   Composite Nonconvex Stochastic Optimization: Block-specific Steplengths and Adapted Batch-sizes}

 \author{Jinlong \ Lei  \and  Uday \ V. \ Shanbhag
 \thanks{ Lei is  the Department of Control Science and Engineering, Tongji University; She is also with the Shanghai Institute of Intelligent Science and Technology, Tongji University, Shanghai 200092, China {\tt \small leijinlong@tongji.edu.cn}.}
 \thanks{ Shanbhag is with
the Department of Industrial and Manufacturing Engineering, Pennsylvania State
University, University Park, PA 16802, USA   {\tt\small udaybag@psu.edu}.} 
\thanks{ The work  has been partly supported  by NSF grant 1538605	and 1246887 (CAREER).}}

\maketitle

\begin{abstract}
This work considers  the minimization of a  sum of an expectation-valued
coordinate-wise $L_i$-smooth  nonconvex  function  and  a nonsmooth
block-separable convex regularizer where $i$ denotes the block index.
\red{We  propose an asynchronous variance-reduced algorithm,} where in each
iteration, a single  block is randomly chosen to update its estimates  by  a
proximal   variable sample-size stochastic gradient scheme,  while the
remaining blocks are kept invariant. Notably, each block employs a steplength
that is in accordance with its block-specific Lipschitz constant while
block-specific batch-sizes are random variables    updated at a rate that grows
either at a geometric or a polynomial rate with the (random) number of times that
block is selected. We show  that every limit point for almost every sample path
is a stationary point and  establish the ergodic non-asymptotic rate of
$\mathcal{O}(1/K) $  in terms of expected sub-optimality.  Iteration and oracle complexity to obtain an
$\epsilon$-stationary point are shown to be $\mathcal{O}(1/\epsilon)$  and
$\mathcal{O}(1/\epsilon^2)$, respectively.  Furthermore, under a $ \mu
$-proximal Polyak-{\L}ojasiewicz  condition with the batch size increasing at a
geometric rate,  we prove that the suboptimality diminishes at a {\em
geometric} rate, the {\em optimal} deterministic rate while iteration and
oracle complexity to obtain an $\epsilon$-optimal solution are proven to be
$\mathcal{O}( (L_{\rm max}/\mu) \ln(1/\epsilon))$ and $\mathcal{O}\left((L_{\rm
ave}/\mu) (1/\epsilon)^{1+c} \right)$, respectively \red{where $L_{\max}$ and $L_{ave}$ denote the maximum and average of the block-specific Lipschitz constants and $c \geq 0$}. In the
single block setting, we obtain the {\em optimal}  oracle complexity
$\mathcal{O} (1/\epsilon)$.   In pursuit of less aggressive sampling rates,
when the batch sizes increase  at a polynomial rate, \red{we show that the
suboptimality decays at a corresponding polynomial rate  and establish   the
iteration and oracle complexity as well.} Finally, preliminary numerics support
our theoretical findings, displaying significant improvements over schemes
where steplengths are based on global Lipschitz constants.  \end{abstract}


\section{Introduction}

In this  paper, we consider the following  composite  stochastic programming:
  \begin{equation}\label{problem1}
\vspace{-0.1in}
\min_{x\in \mathbb{R}^d } \ F(x )\triangleq \bar{f}(x)+\sum_{i=1}^n r_i(x_i) ,
\end{equation}
where   $x_i\in\mathbb{R}^{d_i}$,   $x  \in \mathbb{R}^d $ is  partitioned  into  $n$
blocks as  $x=(x_1,\cdots,x_n)  $  with $d \triangleq \sum_{i=1}^n d_i$,  $r_i:\mathbb{R}^{d_i} \to \mathbb{R}$ is a   convex  nonsmooth  function,  $\bar{f}(x)\triangleq \mathbb{E}_{\xi}[f(x_1,\cdots,x_n,\xi)]$
is an expectation-valued   smooth possibly nonconvex function  with
coordinate-wise  $L_i$-Lipschitz continuous gradients,  the random
vector  $\xi: { \Omega} \to \Real^m$  is defined on the probability
space $({ \Omega}, {\cal F}, \mathbb{P})$, and $f :
\mathbb{R}^d\times \mathbb{R}^m \to \mathbb{R}$ is a scalar-valued
function.  Suppose that     problem  \eqref{problem1} has at least one solution.  Let $X^*$ and  $F^*$ denote the  set of optimal solutions  and the optimal function value, respectively.
 Nonsmoothness  may   be addressed through the   proximal  operator~\cite{rockafellar1976monotone,parikh2014proximal}, defined as
\begin{equation}\label{proximal}
\textrm{prox}_{\alpha r}(x) \ \triangleq \ \argmin_{y} \left( r(y)+{1\over 2\alpha} \|y-x\|^2\right)   ,
\end{equation}
where $r(\cdot)$ is a closed and convex function,   $ \alpha>0,$ and  the  argmin is uniquely defined.
We will propose an asynchronous  proximal  variance-reduced   block scheme for solving  the problem \eqref{problem1}. 

 \red{ \subsection{ Prior research.}}  We first  review the literature on proximal, variance-reduced, and block methods.

  {\em (i) Proximal-gradient  methods.}
 Proximal-gradient  methods and  their
accelerated variants  are among the most important methods for solving composite convex  problem  $ f(x)+r(x)$   (also see
forward-backward splitting (FBS) methods~\cite{lions1979splitting,combettes2011proximal,mine1981minimization}).
 While accelerated (or unaccelerated) schemes \cite{beck2009fast} display non-asymptotic convergence rates in function value of $\mathcal{O}(1/k^2)$ (or $(\mathcal{O}(1/k))$),
FBS methods~\cite{chen1997convergence,tseng2000modified}  display linear convergence when $\nabla f (x)$ is   strongly monotone. Nonconvex  extensions have been studied  in~\cite{attouch2013convergence,frankel2015splitting,karimi2016linear}, where the convergence   to a stationary  point is shown in~\cite{attouch2013convergence} while rate  statements are  provided under both  the Kurdyka-{\L}ojasiewicz
(KL) property~\cite{frankel2015splitting} and the  Polyak-{\L}ojasiewicz (PL)
condition~\cite{karimi2016linear} (where a linear rate is
proven).

  {\em  (ii)  Variance\red{-reduced} schemes.} A
{\bf stochastic   proximal gradient} method  was presented in
\cite{rosasco2014convergence}  for solving  composite convex stochastic
optimization, where the a.s. convergence and   a mean-squared convergence rate
$\mathcal{O}(1/k)$ were developed  in strongly convex regimes,   in sharp
contrast with the linear rate of convergence  in  deterministic settings.
Variance reduction \red{and variable sample-size} schemes  have gained increasing relevance in first-order
methods for stochastic convex
optimization~\cite{shanbhag2015budget,jalilzadeh2018optimal,ghadimi2016accelerated,ghadimi2016mini,
jofre2017variance}; in the class of \red{variable sample-size} schemes, the true gradient
is replaced by the  average of an increasing  batch of sampled gradients,
 progressively reducing the variance of the sample-average.  In strongly convex regimes, linear rates
were  shown for stochastic gradient\ methods~\cite{shanbhag2015budget,jofre2017variance} and extragradient
methods~\cite{jalilzadeh2018optimal}, while  for merely convex optimization problems, accelerated rates of
$\mathcal{O}(1/k^2)$ and $\mathcal{O}(1/k)$ were proven for
smooth~\cite{ghadimi2016accelerated,jofre2017variance} and
nonsmooth~\cite{jalilzadeh2018optimal} regimes, respectively. Mini-batch stochastic approximation   methods were developed in \cite{ghadimi2016mini} for nonconvex stochastic composite optimization.
Alternative  \red{variance reduction}  schemes like SAGA and SVRG,   applied   to  finite-sum machine learning
problems~\cite{reddi2016stochastic,reddi2016fast,xiao2014proximal,reddi2016proximal},
rely on periodic use of the exact gradient, leading to  recovery of
deterministic convergence rates.  For example, geometric rates were provided by~\cite{reddi2016proximal} for
\red{proximal   SAGA and   SVRG  algorithms}   in nonconvex regimes under the proximal PL inequality.
\red{The variance reduced schemes have   been  widely applied to the data-driven  applications, e.g., the
phase retrieval  problem   \cite{wang2017scalable,wang2016solving}  that aims at  reconstructing a general signal vector from magnitude-only measurements. }

  {\em  (iii) Block coordinate descent (BCD) schemes.} BCD  methods~\cite{d1959convex}  have been   widely
used in machine learning   and optimization,  where variables are
partitioned  into manageable blocks and in  each iteration,  a single block is
chosen to update while the remaining blocks  remain fixed.  Recently,
in~\cite{peng2016coordinate},  coordinate-friendly operators were investigated
that perform low-cost coordinate updates and it is shown that a variety of problems in
machine learning  can be efficiently  resolved by such an update.
The convergence  properties of cyclic BCD methods has been extensively analyzed
in~\cite{tseng2001convergence,razaviyayn2013unified,xu2013block}. Nesterov
considered a {\em randomized}  BCD   method~\cite{nesterov2012efficiency} and
proved sublinear and linear convergence in terms of expected objective value
for general convex  and strongly convex cases, respectively.
In~\cite{richtarik2014iteration},    proximal  (but unaccelerated) extensions
were developed to contend with composite problems (also see
\cite{richtarik2014iteration,davis2016asynchronous,xu2013block,xu2017globally,yousefian18stochastic}),
while in \cite{fercoq2015accelerated}, an accelerated proximal
RBCD  scheme was presented with a rate of $\mathcal{O}(1/k^2).$
 \red{In more recent work \cite{csiba2017global},}  diverse block
selection  rules are considered and linear statements are provided for deterministic
nonconvex problems under the PL  condition.
\red{\subsection{  Motivations.} } We consider a class of techniques that combine
variance reduced  and block-based schemes for solving the {\bf nonconvex}
nonsmooth stochastic programs, drawing inspiration from two seminal papers.  Of
these, the first by Xu and Yin~\cite{xu2015block} proposes a block stochastic
gradient (BSG) method  that cyclically updates blocks of variables. The second
paper, by  Dang and Lan~\cite{dang2015stochastic}, presents a stochastic block
mirror-descent  scheme reliant on  randomly choosing and updating a single
block by a  mirror-descent stochastic approximation method.
In~\cite{xu2015block} and  \cite{dang2015stochastic}, rates are provided in the
convex setting while in nonconvex regimes, Dang and
Lan~\cite{dang2015stochastic} present non-asymptotic rates.  Yet, there are
several {\em shortcomings} that motivate the present research: (1)  {\bf Centralized
batch sizes.} The schemes in ~\cite{xu2015block,dang2015stochastic}
require a {\em  centrally specified batch-size}
across all blocks requiring global knowledge of the  global clock, i.e.  iteration $k$; (2)  {\bf Shorter steps.} Block-invariant steplengths utilize the global Lipschitz constant   leading to significantly  shorter steps and poorer performance; (3)  {\bf No almost sure convergence  guarantees.}
Almost sure convergence  guarantees are unavailable for BCD schemes for general nonconvex problems; (4)  {\bf Sub-optimal rate statements.} Optimal deterministic rates via variance-reduced schemes are unavailable but
have been alluded to in convex regimes~\cite[Rem.~7]{xu2015block};
Refinements via the PL condition remain open
questions.

 \begin{table}  [!htb]
 \centering
\newcommand{\tabincell}[2]{\begin{tabular}{@{}#1@{}}#2\end{tabular}}
 \footnotesize
 \centering

	SC,C,NC: Strongly convex, convex, nonconvex; p-PL: proximal P-L, $\delta \geq 0$
	\\
 \begin{tabular}{c|c|c|c|c|}
        \hline
         &  App. &Metric & Asym/Rate/complexity & Comments     \\ \hline
   \multirow{3}{*}{\cite{xu2015block}} &NC &-- & $\mathbb{E}[ d(0,\partial F(x(k)) ] \to 0$ & \\ \cline{2-4}
 & C & $\mathbb{E}[ F(x(k))-F^*]$ &  $\mathcal{O}(1/\sqrt{k})$ &{\tabincell{l}{(i) Centralized batch-size $N_k$;\\ (ii) Central. step depend. on global $L$;\\ (iii) No rate for nonconvex; }}
\\ \cline{2-4}
 & SC & $\mathbb{E}[\|x(k)-x^*\|^2]$ &  $\mathcal{O}(1/k)$ &\tabincell{l}{(iv) sub-linear rates for SC} \\ \hline
    \multirow{3}{*}{  \cite{dang2015stochastic}} &NC & $\mathbb{E} [ \|G_{\alpha}(x_{\alpha,K})\|^2  ] $&   \tabincell{l} {rate: $\mathcal{O}(1/K)$\\
iteration:   $\mathcal{O}\left(  n  L_{\max}/\epsilon  \right) $\\
oracle: $\mathcal{O}\left(n^2  \nu^2   L_{\max} /\epsilon^2\right) $ }  &  {\tabincell{l}{(i) Centralized  batch-size $N_k$ ; \\ (ii) Central stepsize depend  on $L_{\max}$; \\
(iii) Inferior constants dependent  on $L_{\rm \max}$; }}\\   \cline{2-4}
 & C & $\mathbb{E}[ F(\bar{x}_K)-F^*]$ &  $\mathcal{O}(1/\sqrt{K})$,
$\mathcal{O}(1/K)$ (SC)  & \\
\hline \multirow{12}{*}{ \tabincell{l} {This \\ work}  } &   \multirow{5}{*}{ NC}  & --&a.s.  conv.    to stationary points   &  (i)  {\bf a.s. convergence (\uvs{unavailable in~\cite{xu2015block,dang2015stochastic}})}
\\ \cline{3-4}&&    $\mathbb{E} [ \|G_{\alpha}(x_{\alpha,K})\|^2  ] $  &    \tabincell{l} {rate: $\mathcal{O}(1/K)$ \\
iteration:  $\mathcal{O}\left(  n L_{\rm ave}/\epsilon  \right) $  \\
oracle:  $\mathcal{O}\left(n^2  \nu^2 L_{\rm ave} /\epsilon^2\right) $  } & {\tabincell{l}{
(ii) Block-specific stepsize   depend  on $L_i$  \\
$\implies$  {\bf better empirical behavior};  \\  (iii)  {\bf Optimal rates} depend. on $L_{\rm ave}$ \\ (\uvs{\bf Rather than $L$, $L_{\max}$~\cite{xu2015block,dang2015stochastic}}) }}
     \\ \cline{2-4}
 &   \multirow{5}{*}{p-PL}  & \multirow{5}{*}{$\mathbb{E}[ F(x(k))-F^*]$}&
  \tabincell{l} {  {\bf geometric batch-size}: \\ geometric  rate \\ iteration:
$\mathcal{O} \left(   {   n L_{\max} \over \mu  } \ln( 1/\epsilon ) \right)  $
\\ oracle:     $\mathcal{O}\Big({n L_{\rm ave } \over \mu}(1/\epsilon)^{
1+\delta }\Big)$    } &  {\tabincell{l}{   \\
(iv) Block-specific  random.~batch-sizes \\
$\implies$ {\bf no central coordination of   clocks} \\ \\ (v) Optimal geometric rate  under p-PL
  }}  \\ \cline{4-4} &&&
 \tabincell{l} {   {\bf poly.~batch-size: deg. $v$}: \\  polynomial rate $\mathcal{O}(k^{-v})$  \\ iteration:  $\mathcal{O} (  (1/\epsilon)^{1/v})$ \\ oracle: $\mathcal{O} \left(  (1/\epsilon)^{1+1/v}  \right)$  }&   \\ \hline     \end{tabular}
 \vskip 2mm \caption{\small List of  literature on block-structured stochastic optimization} \label{TAB-lit}
\end{table}

\red{\subsection{Summary of  Contributions.}}
We address the \red{aforementioned}   gaps through \red{designing a  novel
algorithm} that combines  a  randomized BCD   method  with a proximal VSSG
method, reliant on {\em block-specific steplengths} based on locally available
$L_i$ without    knowledge  of the central Lipschitz constant,    leading to
{\bf larger steplengths and improved behavior},
and on  {\em random block-specific batch-sizes} (adapted to its block-selection
history)    without  knowledge of the  global clock, leading to {\bf  lower
informational  coordination requirements.} We make the following {\em
contributions} supported by numerics in Section~\ref{Sec:exa}.
Table~\ref{TAB-lit} formalizes the distinctions in our scheme, while
Table~\ref{TAB-comp} compares our results with deterministic rates for
nonconvex regimes.

 \begin{table}
\newcommand{\tabincell}[2]{\begin{tabular}{@{}#1@{}}#2\end{tabular}}
 \footnotesize
\centering
{\center{(a) Iteration complexity in smooth case $(r(x)=0)$ \\
\centering
 \begin{tabular}{|c|c|c|c|  }
 \hline  \multicolumn{2}{|c|}{ block selection rule} &  PL   & general nonconvex
 \\ \hline   \hline
 \multirow{2}{*}{unif.} & deterministic \cite{csiba2017global}& \tabincell{c} { $\mathcal{O}\left({n L_{\max}\over \mu}\ln \left({F(x_1)-F^*\over \epsilon} \right) \right)$   } & $\mathcal{O}\left( {nL_{\max} \left(   F(x(0))   -F^*   \right) \over \epsilon} \right) $
  \\ \cline{2-3}
 &  stoch. (This work) & $\mathcal{O}\left({n L_{\max}\over \mu}\ln \left({F(x_1)-F^*\over \epsilon}+ {n \nu^2\over \epsilon} \right) \right)$&
 $\mathcal{O}\left( {nL_{\max} \left(   F(x(0))   -F^*   \right) \over \epsilon}   \right) $ \\ \hline
 \multirow{2}{*}{non-unif.} & deterministic  \cite{csiba2017global} & \tabincell{c} { $\mathcal{O}\left({n L_{\rm ave}\over \mu}\ln \left({F(x_1)-F^*\over \epsilon} \right) \right)$   } & $\mathcal{O}\left( {nL_{\rm ave} \left(   F(x(0))   -F^*   \right) \over \epsilon} \right) $   \\ \cline{2-3}
 &  stoch. (This work)  &  $\mathcal{O}\left({n L_{\rm ave}\over \mu}\ln \left({F(x_1)-F^*\over \epsilon}+ {n \nu^2\over \epsilon} \right) \right)$&
  $\mathcal{O}\left( {nL_{\rm ave} \left(   F(x(0))   -F^*   \right) \over \epsilon}   \right) $ \\ \hline
 \end{tabular}
}  }
\vskip 2mm
{ \center{ (b) Iteration complexity in nonsmooth case $(r(x)\neq 0)$}
\\
\centering
\begin{tabular}{|c|c|c|c|  }
 \hline   \multicolumn{2}{|c|}{ block selection rule} &  PL   & general nonconvex
 \\ \hline   \hline
 \multirow{2}{*}{unif.}  & deterministic\cite{csiba2017global}  & \tabincell{c} { $\mathcal{O}\left({n L_{\max}\over \mu}\ln \left({F(x_1)-F^*\over \epsilon} \right) \right)$    } & $\mathcal{O}\left( {nL_{\max} \left(   F(x(0))   -F^*   \right) \over \epsilon} \right) $
  \\ \cline{2-3}
 &   stoch. (This work)  &  $\mathcal{O}\left({n L_{\max}\over \mu}\ln \left({F(x_1)-F^*\over \epsilon}+ {n \nu^2\over \epsilon} \right) \right)$&
 $\mathcal{O}\left( {nL_{\max} \left(   F(x(0))   -F^*   \right) \over \epsilon}   \right) $ \\ \hline
 \end{tabular}
 }
\vskip 2mm \caption{\small Comparison  with   deterministic rates for nonconvex block methods }\label{TAB-comp}

\end{table}
\vskip 2mm
\noindent (I)  In Section \ref{Sec-3}, we prove that  every limit point for almost every sample path is a stationary point
under appropriately chosen batch sizes  and  show  that the
ergodic mean-squared error of the  gradient mapping diminishes at a rate    $ \mathcal{O}(1/K)$.  We
then establish that  for any given $\epsilon>0,$ the  iteration complexity (no.
 of proximal  evaluations) and oracle complexity  (no. of sampled
gradients)  to obtain an $\epsilon$-stationary point   are $\mathcal{O}(n
L_{\max}/\epsilon ) $  and   $\mathcal{O}(n^2  \red{\sigma}^2 L_{\max}^2 L_{\min}^{-1}  \epsilon^{-2})$   with    uniform  block  selection, where $L_{\max} \triangleq \max_i L_i $,   $L_{\min}
\triangleq \min_i L_i $, \red{and $\sigma^2$ denotes a uniform bound on the variance of gradient noise.}   When the   blocks are  chosen as per a
non-uniform distribution with probabilities $L_i
\left(\sum_{i=1}^nL_i\right)^{-1}$ for each $i=1,\cdots,n$, the iteration   and
oracle complexity  are  improved to $\mathcal{O}(n L_{\rm ave}/\epsilon ) $  and
$\mathcal{O}(n^2 \red{\sigma}^2 L_{\rm ave}/\epsilon^2 ) $ with $ L_{\rm ave} \triangleq
\sum_{i=1}^n L_i/n$.  This represents a constant factor improvement in the rate
from $L_{\max}$ (in~\cite{dang2015stochastic}) to $L_{\rm ave}$. Preliminary numerics reveal that block-specific steplengths lead to significant improvements in empirical behavior (see Table~\ref{tab5}).  \us{In addition, we show that if the batch-sizes increase at a quadratic rate with the number of times a particular block is chosen, the choice of batch-size sequences no longer relies on the knowledge of global Lipschitzian information.}

\noindent (II)   In Section \ref{Sec-4},  we consider a class  of nonconvex
functions satisfying the {\bf  proximal PL condition}  with parameter $\mu$
(see Assumption \ref{ass-PL}) and prove  that  when the  block-specific batch size is
random and  increases at a suitable  geometric rate  with the number of times
the block is selected, the expectation-valued optimality gap
  diminishes  at a {\em geometric} rate. In addition,
with    uniform block  selection,  the iteration and oracle complexity to
obtain an $\epsilon$-optimal  solution  are $ \mathcal{O} \left( {   (n
L_{\max} \slash \mu)  } \ln( 1/\epsilon ) \right) $  and $\mathcal{O}\Big({n
L_{\rm ave } \over \mu} (1/\epsilon)^{  \big(1+{1\over n\kappa_{\min}-1}  \big)
{L_{\max}  \over L_{\min} }}\Big)$    respectively, where
$\kappa_{\min}\triangleq  L_{\min} /\mu$.  While in the smooth  regimes with a
non-uniform  block   selection,   the  iteration and oracle complexity  bounds
are  improved to     $\mathcal{O}\left({n L_{\rm ave}\over \mu}\ln
\left({1\over \epsilon} \right) \right)$  and $\mathcal{O}\Big({n L_{\rm ave }
\over \mu}(1/\epsilon)^{  \big(1+{1\over n\kappa_{\min}-1}  \big) {L_{\rm ave}
\over L_{\min} }}\Big)$,  respectively.  Specifically, when  $n=1$,   the {\bf
optimal} oracle complexity $\mathcal{O} \left({ L \over \mu \epsilon} \right)
$ is obtained.  Notably, these rates match the deterministic versions
in~\cite{csiba2017global}.  \red{We further show that   when the batch size increases
at a    polynomial rate of degree $v\geq 1$, the convergence rate is   $ \mathcal{O}(k^{-v})$ and   the corresponding
iteration and oracle complexity    are  respectively $\mathcal{O} ( v(1/\epsilon)^{1/v})$ and   $\mathcal{O}
\left(e^v v^{2v+1}(1/\epsilon)^{1+1/v}  \right)$.}

   \section{Asynchronous  Block   Proximal Stochastic Gradient  Algorithm} \label{sec:randomized}
We assume access to a {\em proximal oracle} (PO) that  outputs
$\textrm{prox}_{\alpha r_i}(x_i)$ at any point $x_i \in \mathbb{R}^{d_i}$ for any
$\alpha>0.$ Since the exact gradient $\nabla \bar{f}(x)$ is unavailable in a
closed form,  we assume  there exists a {\em stochastic first-order oracle}
(SFO) such that  for every  $i \in \mathcal{N}\triangleq\{1,\cdots,n\}$ and for any given $x,\xi$, a
sampled gradient $\nabla_{x_i} f(x,\xi)$ is returned, which   is  an unbiased
estimator of  $\nabla_{x_i} \bar{f}(x)$. We aim to develop  efficient
algorithms for obtaining an $\epsilon$-optimal solution, where the efficiency
is measured  by the  iteration complexity (no. of PO  calls) and  the oracle
complexity (no. of SFO calls).  Time is slotted at $k=0,1,2,\dots$.  Block $i$
at time $k$ holds a state $x_i(k) \in\mathbb{R}^{d_i}$ that is  an estimate for
the corresponding coordinates of the optimal solution.  We propose an asynchronous  variance-reduced block  stochastic gradient scheme (Algorithm \ref{alg-prox-gradient}) where at time instant $k,$
a block $i  \in \mathcal{N}$ is randomly  chosen with probability $p_{i }$ to
compute  the proxima  update    \eqref{proxg}, where $\alpha_i$ is the  constant steplength and $N_i(k)$  is the number  of sampled gradients utilized at block $i$, respectively.   To be specific, the steplength  $\alpha_i$  depends on   its block-specific Lipschitz constant $L_i$    and  the batch-size $N_i(k)$ is a function of the  random number of times   block $i$ is selected up to time $k.$  We will specify the selections of   $\alpha_i$ and $N_i(k)$ upon the performance analysis in Sections \ref{Sec-3} and \ref{Sec-4}.
\begin{algorithm}
\caption{ Asynchronous  variance-reduced block  stochastic gradient algorithm }\label{alg-prox-gradient}
 Let $k:=0$, $ x_i(0) \in  \mathbb{R}^{d_i} $    and $0 < p_i < 1$ for $i = 1, \hdots, n$ such that
$\sum_{i=1}^n  p_i = 1.$
\begin{enumerate}
\item[(S.1)] Pick $i_k=i \in \mathcal{N}$ with probability $p_i$.
\item[(S.2)] If $i_k=i$, then  block  $i$  updates  its state  $ x_i(k+1) $  as follows:
\begin{equation}\label{proxg}
x_i(k+1) =\textrm{prox}_{ \alpha_i r_i} \left(x_i(k)- \alpha_i { \sum_{j=1}^{N_i(k)}\nabla_{x_i} f(x(k),\xi_j(k)) \over N_i(k)} \right),
\end{equation}
 where   $\alpha_i>0$ is the steplength of block $i$,   $N_i(k)$ is   the number  of sampled gradients utilized at block $i$,  and samples $\{ \xi_j(k)\}_{j=1}^{N_i(k)}$ are  randomly generated  from the  probability space $({ \Omega}, {\cal F}, \mathbb{P})$.  Otherwise,     $ x_j(k+1):= x_j(k)$ if $j  \neq  i_k$.
 \item[(S.3)] If $k > K$, stop and return $\{x(k)\}_{k=0}^K$; Else, $k:=k+1$ and return to (S.1).
\end{enumerate}
\end{algorithm}

\begin{remark}
 In fact, Algorithm  \ref{alg-prox-gradient} does  not require a global coordinator  to coordinate the  block selection.   The  block selection rule   (S.1.)   accommodates the  Poisson model  employed by \cite{boyd2006randomized}  as   a special case,  where   each block $i$  is activated  according to   a local Poisson clock that  ticks according to a Poisson process with rate $\varrho_i>0$. Suppose that the local Poisson clocks are independent  and  there is a virtual global clock which ticks whenever
any of the local Poisson clocks tick, then  the global clock ticks
according to a Poisson process with rate $\sum_{i=1}^n\varrho_i$.
Let  $Z_k$ denote the time of the $k$-th tick of the global
clock.  Since the  local Poisson clocks are independent,  with
probability one,   there is  a single block  whose Poisson clock  ticks at time $Z_k
$  with probability  $\mathbb{P}(i_k =i)={\varrho_i \over
\sum_{i=1}^n\varrho_i} \triangleq p_i.$     In addition, note by (S.2.) and the  algorithmic  parameter selection rule  that each block  $i$  maintains  a local clock counting the  number of its block updates, and that the update   \eqref{proxg}
does not necessitate knowledge of either the global Lipschitz constant  or the  virtual global clock.
 As such,  Algorithm \ref{alg-prox-gradient} is an  asynchronous scheme with limited coordination across  blocks.
Thus,  in practical  applications, Algorithm  \ref{alg-prox-gradient} might be   helpful  for the  decentralized implementations  compared to those designed in  \cite{dang2015stochastic} and \cite{xu2015block}.

\end{remark}

If the observation noise $w_i(k)$ of the exact gradient is defined as 
  \begin{align}\label{grad_noise} w_i(k+1)\triangleq  \tfrac{\sum_{j=1}^{N_i(k)}\nabla_{x_{i }} f(x(k),\xi_j(k)\red{)}}{N_i(k)}-\nabla_{x_{i }} \bar{f}(x(k)),\end{align}
 then \eqref{proxg} may be rewritten as
  \begin{equation}\label{proxg-rw}
x_i(k+1) =\textrm{prox}_{\alpha_i r_i} \big(x_i(k)-\alpha_i \left(\nabla_{x_{i}} \bar{f}(x(k))+w_i(k+1)\right) \big).
\end{equation}
By  taking   $r_i(x_i) $  as an indicator function of a  convex set $X_i  $, i.e., $r_i(x_i) = 0 $ if $ x_i\in X_i $ and $r_i(x_i) =+\infty  $ otherwise,     then    the problem \eqref{problem1} reduces to the   stochastic  programming $
\min_{x_i\in X_i, i\in \mathcal{N} }  \mathbb{E}\left[f(x_1,\cdots,x_n,\xi)\right] .$
In this case,     the update \eqref{proxg} reduces to the  variable sample-size projected stochastic gradient method:
\begin{equation*}
x_i(k+1) = \mathbf{P}_{X_i} \left(x_i(k)- \alpha_i \left(\nabla_{x_{i}} \bar{f}(x(k))+w_i(k+1)\right)  \right),
\end{equation*}
where $ \mathbf{P}_{X_i}(x_i) $ denotes the  projection of $x_i$ onto the set $X_i.$
This  can be thought  as a  generalization of the schemes
   proposed  in \cite{shanbhag2015budget,jalilzadeh2018optimal} for  solving  constrained  stochastic convex   program.

\red{We denote the $\sigma$-field of the entire information used by Algorithm \ref{alg-prox-gradient} up to (and
including) the update of $x(k)$ by $\mathcal{F}_k.$  Then $N_i(k)$ is adapted to  $\mathcal{F}_k$ when  $N_i(k)$ is a function of the  random number of times   block $i$ is selected up to time $k.$}
We impose the following  conditions on the  objective functions and  observation noises.

  \begin{assumption}\label{ass-fun}
   (i) $r_i$ is a  proper lower semicontinuous and convex function with effective domain
  $\mathcal{R}_i$ required to be  compact.
(ii) There exists a constant $L_i>0$ such that for any $x_i'\in  \mathcal{R}_i $ and any  $x \in \Pi_{j=1}^n \mathcal{R}_j $, $ \| \nabla_{x_{i}}  \bar{f}(x )- \nabla_{x_{i}}  \bar{f}(x_1,\dots,x_{i-1}, x_i', x_{i+1},\dots, x_n) \|  \leq L_i \| x_i-x_i'\|.$
  \end{assumption}

  \begin{assumption}\label{ass-noise}
   (i)  There exists \red{$\sigma >0$} such that  for any $  i \in \mathcal{N} $ and all $k\geq 1$,  $\mathbb{E}[\|w_i(k+1)\|^2| \mathcal{F}_k]\leq {\red{\sigma^2} \over N_i(k)} ~a.s. $;
(ii) $i_k$ is independent of $\mathcal{F}_k$   for all $k\geq 1.$
  \end{assumption}

\begin{remark}
Assumption \ref{ass-fun}(i) requires the effective domain of each block's  nonsmooth convex regularizer  to  be a compact set. For example, if  $r_i(x_i) $ is set   to be  an indicator function of a  convex set $X_i  $, then it is required  that $X_i$  be  a compact set. Such a condition guarantees that the iterate  sequence  $\{x(k)\}$ produced  by Algorithm \ref{alg-prox-gradient}  is uniformly bounded. Assumption \ref{ass-fun}(ii)  necessitates  the  gradient of the nonconvex smooth  term $\bar{f}(x )$ to be block-wise Lipschitz continuous.   Assumption \ref{ass-noise}(i)  can be satisfied  if  the sampled  gradients are independently   generated  and $\mathbb{E} [\| \nabla_{x_i} f(x,\xi)-\nabla_{x_i} \bar{f}(x )\|^2]$, the variance of the stochastic gradient noise,   is uniformly bounded.
\end{remark}

  \section{Convergence to Stationary Points}\label{Sec-3}
In this section, we will   prove the almost sure convergence of iterates   to a  stationary point    and  establish the non-asymptotic rate of  Algorithm \ref{alg-prox-gradient}.

\subsection{Preliminary Lemmas}
Before presenting the convergence results, we recall a preliminary result from   \cite[Lemma 2]{reddi2016proximal}.
  \begin{lemma}\label{pre-lem}
  Suppose
  $y \triangleq \textrm{prox}_{\alpha r}(x-\alpha g)$
  for some $g\in \mathbb{R}^d.$ Then for any $z\in \mathbb{R}^d,$ 
  \begin{align*}
\bar{f}(y)+r(y)&\leq \bar{f}(z)+r(z)+(y-z)^T(\nabla \bar{f}(x)-g)+\left(\frac{L}{2}-{1\over 2 \alpha}\right) \| y-x\|^2\\&+
\left(\frac{L}{2}+{1\over 2 \alpha}\right) \| z-x\|^2
- {1\over 2 \alpha} \| y-z\|^2.
\end{align*}
  \end{lemma}

We now  give  a simple    relation on the conditional  expectation of  the function value in the following lemma, for which the proof   can be found in Appendix A.  This  is  an important preliminary  result because the  convergence results to be presented   are essentially obtained through a recursive application of this basic lemma.
\begin{lemma} \label{lem-sec3-2}  Suppose that  Assumptions \ref{ass-fun} and  \ref{ass-noise}(ii)   hold. Let  $\{ x(k)\}$ be generated by Algorithm \ref{alg-prox-gradient}, where   $0<\alpha_i\leq{1\over L_i}  $ for each $i=1,\cdots,n$. Define   \begin{equation}\label{def-xbar}
\bar{x}_i(k+1) \triangleq \textrm{prox}_{\alpha_i r_{i}} \left(x_i(k)-\alpha_i \nabla_{x_{i}} \bar{f}(x(k))  \right)
\end{equation}  as the   update if   the true gradient is used.  Then the following recursion holds \red{almost surely (a.s.)}:
 \begin{equation}\label{lem-x-bd}
\begin{split}
\mathbb{E}\big[ F(x(k+1)) |\mathcal{F}_k \big]&\leq  F(x(k))  -\sum_{i=1}^n p_i  \left({1\over 2 \alpha_i}-L_i\right) \|\bar{x}_i(k+1)-x_i(k)\|^2  \\&+   {1 \over 2}     \mathbb{E}\left[   \alpha_{i_k}    \|w_{i_k}(k+1)||^2|\mathcal{F}_k\right], \quad \forall k\geq 1.
\end{split}
\end{equation}
\end{lemma}

Throughout the paper, all inequalities and equalities between random variables
    are assumed to hold a.s., but we often omit to write ``a.s." for simplicity.

\subsection{Asymptotic Convergence}\label{3-sub1}

 {From Assumption \ref{ass-fun}(i) it is seen that the sequence of estimates $\{x(k)\}$ produced  Algorithm \ref{alg-prox-gradient}  is   bounded. We now establish the  almost sure  convergence  by showing that
for  every   limit point of  almost  every   sample path 
$\{x(k)\}$ 
is a stationary point of problem \eqref{problem1}.

\begin{theorem} \label{thm1} Let  $\{ x(k)\}$ be generated by
 Algorithm \ref{alg-prox-gradient}.   Suppose  Assumptions \ref{ass-fun} and  \ref{ass-noise}   hold, and that for each $i=1,\cdots,n$,  $0<\alpha_i<{1\over 2L_i}$ and $\sum_{k=0}^{\infty} {1\over N_i(k)}<\infty ~a.s.$ .
Then   every    cluster point of  almost every sample path $\{x(k)\}$ is a stationary point.
\end{theorem}
{\bf Proof.}
  Note  by Assumption  \ref{ass-noise}(i)  that
 $$ \mathbb{E}\left[   \alpha_{i_k}    \|w_{i_k}(k+1)||^2|\mathcal{F}_k\right] \leq
\sum_{i=1}^n \mathbb{E}\left[   \alpha_i    \|w_i(k+1)||^2|\mathcal{F}_k\right] \leq  \sum_{i=1}^n {\alpha_i \sigma^2\over N_i(k)}.$$
Note that $ {1\over 2 \alpha_i}-L_i>0 $ by   $0<\alpha_i<{1\over 2L_i}$. Then by recalling that $\sum_{k=0}^{\infty} {1\over N_i(k)}<\infty~ a.s.$, we may then apply \cite[Thm.~1]{robbins1985convergence} to  inequality~\eqref{lem-x-bd} and conclude that  $F(x(k)) $ converges almost surely  and $\sum\limits_{k=1}^\infty     \sum_{i=1}^n p_i \left( {1\over 2 \alpha_i}-L_i \right) \|\bar{x}_i(k+1)-x_i(k)\|^2<\infty ~ a.s.$ .  This implies that
 \begin{equation}\label{dif-lim}
\sum\limits_{k=0}^\infty   \|\bar{x}(k+1)-x(k)\|^2   <\infty, ~a.s .
\end{equation}
 Then for almost  every     sample path, we have that $\|\bar{x}(k+1)-x(k)\| \to 0$ as $k \to \infty.$
Let $\hat{x}$ be a cluster point of   any such sequence $\{x(k)\}$.
Then there exists a subsequence  $\{x(k_t)\}$
such that $\lim\limits_{t \to \infty} x(k_t)= \hat{x} $
and hence  $\lim\limits_{t \to \infty} \bar{x}(k_{t}+1)= \hat{x}$  by \eqref{dif-lim}.
For any $i=1,\cdots,n$, using the  definitions \eqref{proximal} and \eqref{def-xbar}, we obtain  that
\begin{equation}\label{barx-mini}
\bar{x}_i(k+1) =\argmin_{y\in \mathbb{R}^{d_i}}\Big[\nabla_{x_i}  \bar{f}(x(k))^T(y-x_i(k))+{1\over 2\alpha_i } \|y-x_i(k)\|^2+ r_i(y)  \Big].
\end{equation}
Then by using the first-order optimality condition, we obtain for all $t:$
\begin{equation}\label{fo-opt}
-{1\over \alpha_i }\left(\bar{x}_i( k_{{t}}+1)- {x_i(k_t)}\right) \in \nabla_{x_i} \bar{ f}\left({x(k_t)}\right)
+\partial r_i\left( {\bar{x}_i( k_t+1)}\right)  .
\end{equation}
By passing to  the limit in \eqref{fo-opt},  using
$\|\bar{x}(k_t+1)-x(k_t)\|\to 0~a.s. ,$ by  $\lim\limits_{t \to \infty}
 {x(k_t)}=\lim\limits_{t \to \infty}  {\bar{x}(k_t+1)}= {\hat{x}}$, and
the continuity of $\nabla_{x_i}  \bar{f} {(\cdot)}$ and the closedness of
$\partial r_i$,  we obtain that $ 0 \in \nabla_{x_i} \bar{
f}( \hat{x})+\partial r_i( \hat{x}_i) $ for each $i=1 \dots,n.$
Thus,  $ 0 \in \nabla  \bar{ f}(\hat{x})+\partial r ( \hat{x}),$ implying that $ \hat{x}$  is a stationary point of
\eqref{problem1}.   \hfill $\Box$

\begin{remark}  \label{rem2} (i)
For any $i\in \mathcal{N}$ and $ k\geq 1$,  define $ \Gamma_i(k) \triangleq
\sum_{p=0}^{k-1 }I_{[i_p=i]}$ \red{as the number of updates  block $i$ has
carried out up to time $k$}, where $I_{[a=b]}=1$ if $a=b$, and $I_{[a=b]}=0$
otherwise. Thus, $ \Gamma_i(k)$ is adapted to $\mathcal{F}_k$, and
$\sum_{k=0}^{\infty} {1\over N_i(k)}<\infty $ holds almost surely  by setting
$N_i(k) \triangleq \lceil  ( \Gamma_i(k)+1 )^{1+\delta}  \rceil  $  for some
$\delta>0$. From   \cite[Lemma 7]{koshal2016distributed}, it follows   that
for almost every $\omega \in \Omega $,  there exists  a sufficiently large
$\tilde{k}$ possibly  contingent on the sample path   $\omega $ such that  for
any $ k \geq \tilde{k},$ $\Gamma_i(k) \geq     {k p_i \over 2} ,  i=1,\cdots, n
.$

\noindent (ii) If $\bar{f}(x)$ is convex, then   Theorem \ref{thm1} implies that  $F(x(k))$ converges almost surely to  the optimal function value $F^*$, and  {every} cluster point of   {almost every sample path} $\{x(k)\}$ is a global  minimum   to the problem \eqref{problem1}.
\end{remark}

\subsection{Non-asymptotic Rate}\label{3-sub2}
Recall that for convex optimization, a frequently-used metric is the
sub-optimality  metric $F(x)-F^*$ or the distance to the optimal solution set
$d(x,X^*)$.  However, in  nonconvex optimization,  the iterates might converge
to stationary points that  are not necessarily global minima, and as a
consequence,  the  standard metric cannot be applied.  Thus,  one crucial
problem in analyzing Algorithm \ref{alg-prox-gradient} for nonconvex
problems lies in the selection of the convergence criterion.  In smooth
regimes, it is typical to use $\|\nabla \bar{f}(x)\|$  while in nonsmooth
settings,   an appropriate alternative   is   the proximal gradient mapping
\cite{reddi2016proximal}: $
G_{\alpha}(x)={1\over \alpha} \Big(x-\textrm{prox}_{\alpha r } \big(x-\alpha \nabla \bar{f}(x) \big)\Big) $.
Then $x^0\in \mathbb{R}^d$  satisfying   $G_{\alpha}(x^0)=0$  is   a {\em stationary point} of \eqref{problem1}.
 We now analyze the   rate of convergence of Algorithm \ref{alg-prox-gradient},  and establish     iteration and  oracle complexity bounds to  obtain an  $\epsilon$-stationary point, by using  the following metric to measure   stationarity.
\begin{align}\label{metric}G_{i,\alpha_i}(x)={\frac{  x_i- {\rm prox}_{\alpha_i r_{i}} \left(x_{i }-\alpha_i \nabla_{x_{i}} \bar{f}(x )  \right) }{\alpha_i}} , ~\forall i\in \mathcal{N}, ~ \mathbf{G}_{\bm{\alpha}}(x) \triangleq  \left(G_{i,\alpha_i}(x)\right)_{i=1}^n. \end{align}
It is seen that any zero  of  $\mathbf{G}_{\bm{\alpha}}(x) $ is a stationary point of \eqref{problem1}. Next, we establish a result for    Algorithm \ref{alg-prox-gradient} when  the block is chosen according to a uniform distribution.
\begin{theorem} \label{thm-non-asy} Suppose Assumptions \ref{ass-fun} and \ref{ass-noise}   hold.  Let  $\{ x(k)\}$ be generated by Algorithm \ref{alg-prox-gradient}, where  $p_i={1\over n}$ and $\alpha_i={1\over 4 L_i}$ for each  $i=1,\cdots,n $.   Let $x_{\alpha,K}$  be chosen  from $\{x(k)\}_{k=0}^K$ as per a uniform distribution. We  have the following bound on the  mean-squared error:
\begin{align}  \label{x-expect-bd-1}
   & \mathbb{E}\left[   \|\mathbf{G}_{\bm{\alpha}}(x_{\alpha,K})\|^2 \right]      \leq  { 16 n L_{\max} \left(   \mathbb{E}[F(x(0))]   -F^*   \right) \over K+1} +  {2\sigma^2L_{\max}    \over K+1}  \sum_{i=1}^n     \sum_{k=0}^K  L_i^{-1}\mathbb{E} [    N_i(k)^{-1} ] .
\end{align}
 Then for any given $\epsilon >0,$  by setting  $K=\bar{K}_1(\epsilon) \triangleq     {32 n L_{\max} \left(   \mathbb{E}[F(x(0))]   -F^*   \right)\over \epsilon}  $ and $N_i(k)\equiv \bar{N}_1(\epsilon)\triangleq
{4n\sigma^2   L_{\max} \over \epsilon  L_{\min}},$   the iteration and  oracle complexity to obtain an $\epsilon$-stationary point such that $\mathbb{E}\left[ \|\mathbf{G}_{\bm{\alpha}}(x_{\alpha,K})\|^2 \right] \leq \epsilon$ are  $\bar{K}_1(\epsilon)$ and $\bar{K}_1 (\epsilon) \bar{N}_1 (\epsilon) $, respectively.
\end{theorem}
{\bf Proof.}   Note  by  $\mathbb{P}(i_k=i)=p_i$ and Assumption  \ref{ass-noise}(i), we obtain  that
$$ \mathbb{E}\left[   \alpha_{i_k}    \|w_{i_k}(k+1)||^2 \right] =
\sum_{i=1}^n p_i \alpha_i  \mathbb{E}\left[     \|w_i(k+1)||^2 \right] \leq  \sum_{i=1}^n {\alpha_i  p_i \sigma^2 }\mathbb{E}\left[    N_i(k)^{-1}\right].$$
By  taking  unconditional expectations  of \eqref{lem-x-bd} and rearranging the terms,  the following  holds.
 \begin{equation}\label{lem-x-bd2}
\begin{split}
&\mathbb{E} \left[ \sum_{i=1}^n p_i  \left({1\over 2 \alpha_i}-L_i\right) \|\bar{x}_i(k+1)-x_i(k)\|^2 \right]  \\&\leq  \mathbb{E}\big[ F(x(k)) \big] - \mathbb{E}\big[ F(x(k+1)) \big]+  {  1 \over 2}  \sum_{i=1}^n {\alpha_i  p_i \sigma^2 }\mathbb{E}\left[    N_i(k)^{-1}\right].
\end{split}
\end{equation}
Thus, by summing up  \eqref{lem-x-bd2} from $k=0$ to $K$, we have that
 \begin{equation}\label{x-expect-bd}
\begin{split}
  & \sum_{k=0}^K \mathbb{E} \left[ \sum_{i=1}^n p_i  \left({1\over 2 \alpha_i}-L_i\right) \|\bar{x}_i(k+1)-x_i(k)\|^2 \right]\\&  \leq     \mathbb{E}[F(x(0))] -\mathbb{E}\left[ F(x(K+1))   \right] +  { \sigma^2  \over 2}\sum_{k=0}^K\sum_{i=1}^n \alpha_i p_i\mathbb{E}\left[    N_i(k)^{-1}\right] .
\end{split}
\end{equation}
 By the definitions   \eqref{def-xbar} and \eqref{metric}, $\alpha_i={1\over 4 L_i}$,and $p_i={1\over n}$, we have that
\begin{align*}
    & p_i \left({1\over 2 \alpha_i}-L_i\right) \|\bar{x}_i(k+1)-x_i(k)\|^2 =  p_i  \alpha_i  \left({1\over 2 }-\alpha_i L_i\right) \left({\|\bar{x}_i(k+1)-x_i(k)\|\over  \alpha_i }\right)^2
\\&=  {1\over 16 n L_i} \left\| G_{i,\alpha_i}(x(k))  \right\|^2 \geq  {1\over 16 n L_{\max}} \left\| G_{i,\alpha_i}(x(k))  \right\|^2.\end{align*}
This combined with \eqref{x-expect-bd}, $F(x_{K+1}) \geq F^*$,  and $p_i={1\over n}$  implies that
\begin{align*}
{1\over 16 n L_{\max}} \sum_{k=0}^K  \mathbb{E} \left[ \|\mathbf{G}_{\bm{\alpha}}(x(k))\|^2\right]  \leq      \mathbb{E}[F(x(0))]  -F^*  +  {   \sigma^2 \over 8n}  \sum_{i=1}^n{1\over L_i} \sum_{k=0}^K  \mathbb{E}\left[    N_i(k)^{-1}\right] .
\end{align*}
Therefore, by multiplying  both sides  of the above equation by  ${16 n L_{\max} \over K+1}  $,    {by  using
$$ \mathbb{E}\left[ \|\mathbf{G}_{\bm{\alpha}}(x_{\alpha,K})\|^2 \right] = {  1 \over K +1} \sum_{k=0}^K\mathbb{E}\left[ \|\mathbf{G}_{\bm{\alpha}}(x(k))\|^2 \right] $$ we obtain  \eqref{x-expect-bd-1}}. Since $N_i(k)= \bar{N}_1(\epsilon)$,  we have ${1\over K+1} \sum_{k=0}^K   \mathbb{E}\left[    N_i(k)^{-1}\right] = {1 \over \bar{N}_1(\epsilon)}.$ This  combined with  \eqref{x-expect-bd-1} produces the following bound.
\begin{equation}\label{cor-r1}
\begin{split}
\mathbb{E}\left[   \|\mathbf{G}_{\bm{\alpha}}(x_{\alpha,K})\|^2 \right]  \leq  { 16 n L_{\max} \left(    \mathbb{E}[F(x(0))]  -F^*   \right) \over  \bar{K}_1(\epsilon) +1} +    {2n\sigma^2  L_{\max}    \over \bar{N}_1(\epsilon) L_{\min}}   \leq    \epsilon.
\end{split}
\end{equation}
 Since a single block is chosen to  be  updated  at  each   iteration, the total number of samples used to update  $x(K)$
 is $\sum_{k=0}^{K-1} \sum_{i=1}^n N_i(k) I_{[i_k=i]} =  \bar{K}_1 (\epsilon) \bar{N}_1 (\epsilon) . $
Then  the number of  PO and SFO calls  required to ensure that
$\mathbb{E}\left[ \|\mathbf{G}_{\bm{\alpha}}(x_{\alpha,K})\|^2 \right] \leq \epsilon$ are  $\bar{K}_1(\epsilon)$ and $\bar{K}_1 (\epsilon) \bar{N}_1 (\epsilon)  $, respectively.
\hfill $\Box$

\red{In the following, we} analyze  the rate of convergence of Algorithm \ref{alg-prox-gradient}  with  the active block chosen  via a  non-uniform distribution  constructed using block-specific Lipschitz constants.
 \begin{theorem} \label{thm-non-asy2}Suppose  Assumptions \ref{ass-fun} and  \ref{ass-noise}   hold.
 Let  $\{ x(k)\}$ be generated by Algorithm \ref{alg-prox-gradient}, where $p_i={L_i\over \sum_{i=1}^n L_i} $ and  $\alpha_i = {1\over 4L_i}$  for each $i = 1, \hdots,n$. Then
\begin{align}\label{G-xa}
  \mathbb{E}\left[   \|\mathbf{G}_{\bm{\alpha}}(x_{\alpha,K})\|^2 \right]   \leq  { 16 n L_{\rm ave} \left(  \mathbb{E}[ F(x(0))]   -F^*   \right) \over K+1} + {  2\sigma^2 \over K+1}     \sum_{i=1}^n    \sum_{k=0}^K  \mathbb{E}\left[    N_i(k)^{-1}\right].
\end{align}
Thus,  for any given $\epsilon >0,$ by setting  $K=\bar{K}_2(\epsilon)=   {32n L_{\rm ave} \left(  \mathbb{E}[ F(x(0))] -F^*   \right)  \over  \epsilon} $ and $N_i(k)\equiv \bar{N}_2(\epsilon)={4n\sigma^2 \over \epsilon} $,  the iteration and  oracle complexity to obtain an $\epsilon$-stationary point such that $\mathbb{E} [ \|\mathbf{G}_{\bm{\alpha}}(x_{\alpha,K})\|^2 ] \leq \epsilon$ are  $\bar{K}_2(\epsilon)$ and $\bar{K}_2 (\epsilon) \bar{N}_2 (\epsilon) $, respectively.
\end{theorem}
{\bf Proof.} By definitions  \eqref{def-xbar}, $\alpha_i={1\over 4L_i},p_i={L_i\over  n L_{\rm ave } }$, we have that
\begin{align*}
& p_i  \left({1\over 2 \alpha_i}-L_i\right) \|\bar{x}_i(k+1)-x_i(k)\|^2
\\& =p_i \alpha_i  \left({1\over 2 }-\alpha_i L_i\right) \left({\|\bar{x}_i(k+1)-x_i(k)\|\over  \alpha_i }\right)^2
 =   {1\over 16 n L_{\rm ave}}\left\| G_{i,\alpha_i}(x(k))  \right\|^2.\end{align*}
Then using \eqref{x-expect-bd},   the definition  of $\mathbf{G}_{\bm{\alpha}}(x) $ and $\alpha_i p_i={1\over 4 L_{\rm ave}}$, we have the following:
 \begin{equation*}
\begin{split}    {\sum_{k=0}^K  \mathbb{E} \left[ \|\mathbf{G}_{\bm{\alpha}}(x(k))\|^2\right]  \over 16 n L_{\rm ave}}  \leq    \mathbb{E}[ F(x(0))- F(x(K+1))  ] +  { \sigma^2    \over 8 nL_{\rm ave}}\sum_{i=1}^n \sum_{k=0}^K\mathbb{E}\left[    N_i(k)^{-1}\right] .
\end{split}
\end{equation*}
Then by  using $F(x_{K+1}) \geq F^*$  and  multiplying both sides  of the above equation with ${16n L_{\rm ave}\over K+1} $,   we  obtain \eqref{G-xa}. The rest of the proof is the same as that of Theorem \ref{thm-non-asy}.
\hfill $\Box$

 \begin{remark}\label{rem3} We  observe the following regarding Theorems~\ref{thm-non-asy}--\ref{thm-non-asy2}: \\ \noindent (i) Note that   $\mathbf{G}_{\bm{\alpha}}(x)=  \nabla \bar{f}(x)  $ when $r(x)\equiv 0$.
Suppose $\sum_{k=0}^{K} \mathbb{E}\left[    N_i(k)^{-1}\right]$ is bounded for each $i=1,\cdots,n$.
Then we  attain the non-asymptotic rate $\mathbb{E}[\| \nabla f(x_{\alpha,K})  \|^2]=\mathcal{O}(1/K)$,
which is the best known rate for  first-order methods for  deterministic  nonconvex programs~\cite{ghadimi2016accelerated}.

\noindent (ii)  Note  that     iteration and oracle complexity bounds  of Algorithm \ref{alg-prox-gradient}  with uniform  block selection  are respectively $\mathcal{O}(n L_{\max}/\epsilon ) $  and   $\mathcal{O}(n^2  \sigma^2 L_{\max}^2  /(L_{\min}\epsilon^2 ) )$,      while  if the blocks are selected with a likelihood  proportional to the block-specific Lipschitz constant, the bounds  are respectively  reduced to $\mathcal{O}(n L_{\rm ave}/\epsilon ) $  and  $\mathcal{O}(n^2  \sigma^2  L_{\rm ave}/\epsilon^2 )  $.
\red{Therefore, the block-specific selection rule improves the constant in both the iteration and oracle complexity.}

\noindent (iii)      The iteration complexity (no.  of
partial  proximal evaluations)  is   $\mathcal{O}(n/\epsilon ) $. Since the variable is partitioned into $n$ blocks,
 the iteration complexity (no. full proximal evaluations)   is  $\mathcal{O}(1/\epsilon ) $, which is  optimal  for   deterministic gradient  descent methods.

\end{remark}

\blue{In  Theorems  \ref{thm-non-asy} and \ref{thm-non-asy2}, the parameters    $N_i(k)$ and $p_i$ still require some sort of global information  regarding  the Lipschitz constants (namely $L_{\max}, L_{\min}, L_{\rm ave}$) to attain the optimal iteration complexity  $\mathcal{O}(1/\epsilon ) $ and oracle complexity $\mathcal{O}(1/\epsilon^2 )  $.  Such a  parameter selection works effectively  for the setting where   the
prior information about the Lipschitz constants is known to all blocks. 
  To address the situation  where the  batch-size   depends  on neither  $L_{\max}$
 nor $ L_{\min}$,  we add an extra proposition to explore the convergence properties of  Algorithm \ref{alg-prox-gradient} when   the batch-size  $ N_i(k) $ increases according to  a  quadratic function of the block update time $\Gamma_i(k)$. }

\blue{\begin{proposition} \label{prp-nonconvex} Suppose Assumptions \ref{ass-fun} and \ref{ass-noise}   hold. Let  $\{ x(k)\}$ be generated by Algorithm \ref{alg-prox-gradient},
where $N_i(k)=   (\Gamma_i(k) +1) (\Gamma_i(k) +2)  $,  $p_i={1\over n}$ and $\alpha_i={1\over 4 L_i}$ for each  $i=1,\cdots,n $.
 Then for any given $\epsilon >0,$     the iteration and  oracle complexity required to obtain an $\epsilon$-stationary point such that $\mathbb{E}\left[ \|\mathbf{G}_{\bm{\alpha}}(x_{\alpha,K})\|^2 \right] \leq \epsilon$ are respectively $\mathcal{O}(1/\epsilon ) $  and   $\mathcal{O}(1/\epsilon^3  )$.
\end{proposition}}
\noindent {\bf Proof.} \blue{By  definition,   $ \Gamma_i(k) \triangleq \sum_{p=0}^{k-1 }I_{[i_p=i]}  $
 and  $\mathbb{P}(i_p =i)=  p_i$.  Consequently, we have   that
\begin{align}\label{dis-ga}\mathbb{P}(\Gamma_i(k)=m)= \left({k \atop m }\right) p_i^{m } (1-p_i)^{k-m}.\end{align}
It is noticed from  \eqref{dis-ga}    that for all $  i \in \cal{N}:$
 \begin{align}\label{bd-batch}
   & \quad \mathbb{E} \left[\big( \Pi_{t=1}^v (\Gamma_i(k) +t) \big) ^{-1}  \right]    =\sum_{m=0}^{k} \Pi_{t=1}^v (m+t)^{-1}\mathbb{P}(\Gamma_i(k)=m)
  =\sum_{m=0}^{k}   {k! \ p_i^{m } (1-p_i)^{k-m}\over m! (k-m)!  \Pi_{t=1}^v (m+t)}   \notag
  \\& =\sum_{m=0}^{k}   {k! \ p_i^{m } (1-p_i)^{k-m}\over (m+v)! (k-m)!  }   =\Pi_{t=1}^v (k+t)^{-1} p_i^{-v}\underbrace{\sum_{m=0}^{k}   {(k+v)! \ p_i^{m+v } (1-p_i)^{k-m}\over (m+v)! (k-m)!  }}_{ \ \triangleq \ \mathbb{P}(\Gamma_i(k+v) \ \leq \ k)\ \leq \ 1}     \notag
   \\& \leq \Pi_{t=1}^v (k+t)^{-1} p_i^{-v}\leq  (k+1)^{-v} p_i^{-v} .
\end{align}
 Since $N_i(k) \triangleq \prod_{t=1}^2 (\Gamma_i(k)+t)$, it follows that by setting $v = 2$ and by recalling that $p_i^{-v}= n^2$, we have that  $\mathbb{E} [N_i(k)^{-1}]\leq   n^2 \big(  (k+1)(k+2)\big)^{-1}    $.  Incorporating this bound  with \eqref{x-expect-bd-1}
and  noting that  $  \sum_{k=0}^K  \big(   (k+1)(k+2)\big)^{-1} =1-{1\over K +2}\leq 1$ leads to
\begin{align*}
   & \mathbb{E}\left[   \|\mathbf{G}_{\bm{\alpha}}(x_{\alpha,K})\|^2 \right]      \leq  { 16 n L_{\max} \left(   \mathbb{E}[F(x(0))]   -F^*   \right) \over K+1} +  {2 n^2\sigma^2L_{\max}  \sum_{i=1}^n L_i^{-1}  \over K+1}   .
\end{align*}
Hence, the iteration complexity required for achieving $\mathbb{E}\left[ \|\mathbf{G}_{\bm{\alpha}}(x_{\alpha,K})\|^2 \right] \leq \epsilon$ is $K(\epsilon )=  \mathcal{O}(1/\epsilon ) $.

From \cite[p.154]{papoulis1984probability} it follows  that the $t$-th moment of the  binomial distribution $\Gamma_i(k) $
equals the $t$-th  derivative of $M_i(y)$  at $y=0 $, where  $M_i(y)=(p_i e^y+1-p_i)^k $. Thus, we can   show that
 $  \mathbb{E}[\Gamma_i(k)^t]=\mathcal{O}(k^t) $, and hence  $\mathbb{E} \left[  N_i(k)   \right]     =    \mathbb{E}[\Gamma_i(k)^2] +3 \mathbb{E}[\Gamma_i(k)^2] +2=\mathcal{O}(k^2).$  Then      the  expected number of sampled gradients  required to obtain   an   $\epsilon-$stationary point   is   bounded by
 $$\mathbb{E}\Big[  \sum_{k=0}^{ K (\epsilon) -1} \sum_{i=1}^n N_i(k)I_{[i_k=i]} \Big]
  =\sum_{k=0}^{ K (\epsilon)-1} \sum_{i=1}^n p_i \mathbb{E}\left[ N_i(k) \right]  = \mathcal{O} (K^3 (\epsilon))=\mathcal{O}(1/\epsilon^3  ).  \qquad \qquad \qquad  \Box$$

  Proposition \ref{prp-nonconvex} demonstrates  that
for the  general nonconvex problem,    the batch-size $N_i(k)$ increases at a quadratic rate given by $ (\Gamma_i(k) +1) (\Gamma_i(k) +2)  $, where this function  does not   rely  on any global  information  (including the Lipschitz constant
$L$ or the global clock $k$).  This unfortunately leads to poorer bounds on  oracle complexity.  In fact, such a parameter  choice uses the local sampling costs to substitute the cost of global information.  It can   be similarly  shown that when the batch-size is set as   $N_i(k) \triangleq \lceil  ( \Gamma_i(k)+1 )^{1+\delta}  \rceil  $ for some $\delta>0$, then the oracle complexity  is $\mathcal{O}(1/\epsilon^{2+\delta}  )$.

\begin{remark}This work   does not consider the information delay  for
the simplicity of presentation so as to highlight the influence of  the
block-specific  steplengths and   batch-sizes on the   convergence rate.  There
are some recent  papers that study   the asynchronous  schemes with   delays,
e.g.,   asynchronous   algorithms  for  the generalized Nash equilibrium
problem \cite{yi2019asynchronous} as well as   for the stochastic potential game and
nonconvex optimization \cite{lei2019asynchronous}.  However, neither
\cite{lei2019asynchronous}  nor  \cite{yi2019asynchronous} analyze the
convergence rate.  We believe that  the almost sure convergence might still be
guaranteed  by using the uniformly bounded \us{delays}, but possibly at
the cost of a slower convergence rate;  the simulation results displayed in
Figure \ref{fig_delay} validated this thought.\end{remark}}

\section{Global Linear Convergence  under PL-Inequality} \label{Sec-4}

In this section, we will  prove  the global linear convergence of iterates and derive the complexity bounds when the proposed scheme is applied to a class of nonsmooth nonconvex  composite functions satisfying  the  proximal PL inequality. 
 The PL  inequality $ \|\nabla f(x)\|^2 \geq 2\mu \big(f(x)-\min_x f(x)\big)$ requires  the gradient  norm to grow  faster than a quadratic function when moving  away from the optimal value.
It was first proposed  in  \cite{polyak1963gradient}  that  the global  linear convergence of the gradient descent method can be obtained  under the  PL  condition.  Its  generalization,  called the  proximal PL inequality, was proposed in \cite{karimi2016linear} for  the composite function.   It has been shown  in~\cite{karimi2016linear}  that
several important  classes of  functions satisfy  this  proximal PL condition, e.g.,  (i) $\bar{f}$ is strongly convex;
 (ii)  $\bar{f}$ has the form $\bar{f}(x)=h(Ax)$ for
a strongly convex  function $h$ and a matrix $A$ while $r$ is an indicator function for a polyhedral set; and   (iii) $ F$ is convex and satisfies the quadratic growth property. We impose the following  condition on the problem \eqref{problem1}.

  \begin{assumption} ({\bf $\mu$-PL})\label{ass-PL} There is a $\mu>0$ satisfying
${1\over 2} D_r(x,L_{\max})\geq \mu (F(x)-F^*)$,
 where  $$ D_r(x,L)\triangleq -{2 L} \min_{y}\Big[ \nabla \bar{f}(x)^T(y-x)+{ L \over 2  } \|y-x\|^2+r(y)-r(x) \Big].
\vspace{-0.3in}$$
  \end{assumption}
\subsection{Rate Analysis}
We first  present a preliminary lemma, based on which we  show in Theorem \ref{prp1}  and Proposition \ref{prp-rate-poly}
 that $F(x(k)) $ converges in mean to the optimal value $F^*$ at a geometric rate and a polynomial   rate when the  number of the sampled gradients increases at a geometric rate and a polynomial  rate, respectively. The proof of the Lemma \ref{sec4-lem1}
 can be found in Appendix B.
\begin{lemma}\label{sec4-lem1}  Let  $\{ x(k)\}$ be generated by
 Algorithm \ref{alg-prox-gradient}.    Suppose
  Assumptions  \ref{ass-noise} and  \ref{ass-PL}    hold.  Let $\beta \in({1\over 2},1)$ and   $0<\alpha_i\leq {2\beta-1\over L_i(1+\beta)}$.  Define $\alpha_{\min}=\min_{i\in \mathcal{N}}\alpha_i.$ Then for all $k\geq 1$,
\begin{align}
\mathbb{E}\big[ F(x(k+1)) -F^*\big] & \leq \big(1-{\alpha_{\min}(1-\beta) \mu p_{\min} }  \big)\mathbb{E}\big[ F(x(k))-F^*\big]  +{ \sigma^2 \over 2}  \sum_{i=1}^n\alpha_i p_i\mathbb{E}\left[  N_i(k)^{-1}\right].
\label{x-bd3}
\end{align}
\end{lemma}

We now discuss  the  optimal selection of  parameters $\alpha_i$ and $\beta.$  Define $\rho(\alpha,\beta) \triangleq 1-{\alpha(1-\beta) \mu p_{\min} } .$  Then  by   $0< \alpha_i\leq {2\beta-1\over L_i(1+\beta)}$ and  $\beta\in ({1\over 2}, 1)$, there holds $0<\alpha_i  (1-\beta) \leq {(2\beta-1)(1-\beta)   \over L_i(1+\beta)}.$  We set $\beta$ to be the maximizer of $ {(2\beta-1)(1-\beta)  \over 1+\beta }$,  given by $\beta^*=\sqrt{3}-1$.
 Then by setting $ \alpha_i^*= {2\beta^*-1\over L_i(1+\beta^*)}={2-\sqrt{3} \over L_i}$, we get  $ 0<\rho(\alpha^*_{\min},\beta^*)=1- {(2-\sqrt{3})^2 \mu p_{\min}  \over L_{\max}}  <1 .$
  By setting   $\alpha_i= \alpha_i^*$ and $\beta=\beta^* $ in  Algorithm \ref{alg-prox-gradient}, we  obtain the geometric rate    under the proximal PL condition with the geometrically increasing sample-sizes.

\begin{theorem}\label{prp1} {\bf  (Geometric rate of convergence)}
  Suppose $p_i={1\over n}~\forall i\in \mathcal{N}$,  Assumptions   \ref{ass-noise} and  \ref{ass-PL}    hold.
  Consider  Algorithm   \ref{alg-prox-gradient}, where   $ \alpha_i={2-\sqrt{3}\over  L_i}$
 and  $N_i(k)=\left \lceil (1-q_i)^{ -\Gamma_i(k) } \right \rceil$ for some $q_i\in (0,1)$ with  $ \Gamma_i(k) \triangleq \sum_{p=0}^{k-1 }I_{[i_p=i]}$.
Let $q_{\min} \triangleq \min_{i\in \mathcal{N}} q_i$ and $ \rho^* \triangleq  1- {(2-\sqrt{3})^2 \mu   \over n L_{\max}} .$

\noindent (i)  If $ q_{\min} \neq  {(2-\sqrt{3})^2 \mu   \over    L_{\max}}$, then  for all $k\geq 0:$
\begin{align}\label{prp1-r1}
 \mathbb{E}\big[ F(x(k)) -F^*\big]   \leq & \left(1-{1\over n}\min  \left\{  q_{\min}, {(2-\sqrt{3})^2 \mu   \over    L_{\max}} \right\} \right) ^k   \times \notag \\& \Big(    \mathbb{E}\big[  F(x(0))-F^*\big]
+  {\sigma^2 \sum_{i=1}^n \alpha_i  \over 2\left | q_{\min}-{(2-\sqrt{3})^2 \mu /    L_{\max}}\right | }  \Big)  .
\end{align}

\noindent (ii) If $ q_{\min}= {(2-\sqrt{3})^2 \mu   \over  L_{\max}}$, then  the following holds for any  $ \tilde{\rho}\in \left(  \rho^* ,1\right)  $ and all $k\geq 0$:
\begin{equation}\label{linear-rate}
\begin{split}
 \mathbb{E}\big[ F(x(k)) -F^*\big] & \leq \tilde{\rho}^{k } \left(    \mathbb{E}\big[  F(x(0))-F^*\big]
 +{ \sigma^2  \sum_{i=1}^n \alpha _i   \over 2n  \rho^*  \ln \left( \left(\tilde{\rho}/ \rho^*\right)^e\right)  } \right).
\end{split}
\end{equation}
\end{theorem}
{\bf Proof.}  By using \eqref{dis-ga} and $p_i=1/n$, we obtain that     for any $k \geq 1$ and $i \in \cal{N}$:
    \begin{equation}\label{bd-Nk}
\begin{split}
   &\mathbb{E} \left[  N_i(k)^{-1}  \right]   \leq  \mathbb{E} \left[(1-q_i )^{ \Gamma_i(k)}  \right]   =\sum_{m=0}^{k} (1-q_i)^{ m} \mathbb{P}(\Gamma_i(k)=m)   \\&=\sum_{m=0}^{k} \left({k \atop m}\right) (p_i(1-q_i))^m(1-p_i)^{k-m}   =
  (p_i(1-q_i)+1-p_i)^{k}  =\left(1-p_iq_i \right)^{k} .
\end{split}
\end{equation}
Since   $\alpha_i=\alpha_i^*$ and $p_i={1\over n} $, by setting $\beta=\beta^* $, we obtain that
  $\rho(\alpha^*_{\min} ,\beta^* )= 1- {(2-\sqrt{3})^2 \mu p_{\min}  \over L_{\max}}  =  \rho^*$
  with $\alpha^*_{\min} \triangleq \min_{i\in \mathcal{N}} \alpha_i.$
By combining  \eqref{bd-Nk} with   \eqref{x-bd3}  and by  the definition $q_{\min}= \min_{i\in \mathcal{N}} q_i$,
  we have  $
v_{k+1}  \leq   \rho^*   v_k+   \sum_{i=1}^n  {   \alpha_i  \sigma^2 \over 2n }       \left (1-{ q_{\min} \over n}\right)^{k},$
where $ v_k \triangleq \mathbb{E}\big[ F(x(k))-F^*\big] .$  Then  by defining  $q^*\triangleq  1-{ q_{\min} / n} $,  we obtain  that
\begin{equation} \label{linear-opt-gap}
\begin{split}
v_{k+1} & \leq   \left(\rho^*\right) ^{k+1 }  v_0+    \sum_{m=0}^{k}  \left(\rho^*\right) ^m   \left (q^*\right)^{k-m} { \sum_{i=1}^n  \alpha_i \sigma^2 \over 2n }.
\end{split}
\end{equation}

\noindent (i) Suppose  $q_{\min}>{(2-\sqrt{3})^2 \mu   \over   L_{\max}}$,   we then have
 $ q^*= 1-{q_{\min} \over n}<  \rho^*,$  and hence
$$\sum_{m=0}^{k}  \left(\rho^*\right) ^m   \left (q^*\right)^{k-m}   = \left(\rho^*\right) ^{k} \sum_{m=0}^{k}  \left ( q^*/  \rho^* \right)^{k-m}\leq   \left(\rho^*\right) ^{k  }{ 1 \over 1-q^*/  \rho^* } = \left(\rho^*\right) ^{k+1  }{ 1 \over  \rho^*- q^* }  .$$
 Similarly,  for the case  where  $0<q_{\min}<{(2-\sqrt{3})^2 \mu   \over   L_{\max}}$, we have that  $ q^*>\rho^* $  and    $\sum_{m=0}^{k}  \left(\rho^*\right) ^m   \left (q^*\right)^{k-m}   \leq   \left(q^*\right) ^{k+1 }{ 1 \over   q^*- \rho^* }.$
As a result, by using \eqref{linear-opt-gap} we obtain that
$$\mathbb{E}\big[ F(x_{ k+1  }) -F^*\big]  \leq \max  \left \{ q^*, \rho^* \right\}  ^{k+1} \left(     \mathbb{E}\big[ F(x(0))-F^*\big]
+  {   \sum_{i=1}^n \alpha_i \sigma^2 \over 2 n\left | \rho^*- q^*\right | }   \right),$$
and hence \eqref{prp1-r1} holds by definitions of $\rho^*$ and $q^*.$

\noindent (ii) By  $q_{\min}= {(2-\sqrt{3})^2 \mu   \over   L_{\max}}$, we have that $ q^*=\rho^*.$ Then  by \eqref{linear-opt-gap}, we have that $
v_{k }  \leq   \left(\rho^*\right) ^{k  }  v_0+   k\left(\rho^*\right) ^{k } {  \sum_{i=1}^n \alpha_i \sigma^2 \over 2n\rho^* } .$
{By using  $ k \left(\rho^*\right) ^k  \leq  \tilde{\rho}^k   / \ln \left( (\tilde{\rho}/\rho^*)^e\right)$ (see \cite[Lemma 2]{shanbhag2016inexact}) and $ v_0 \triangleq \mathbb{E}\big[ F(x(0))-F^*\big] $, we obtain the result (ii). }
\hfill $\Box$

In some settings,  the evaluation of   sampled gradients might  be costly. Hence  it  is unreasonable  or  impossible to increase the batch-size  too fast at a geometric rate.  As such,  we  try to explore the convergence rate  of Algorithm \ref{alg-prox-gradient}  when  the   batch-size  is  increased at a slower polynomial rate.  Next, we  will investigate  the rate of convergence of Algorithm \ref{alg-prox-gradient}  with  polynomially increasing sample sizes based on a preliminary result from~\cite{lei18game}. 

 \begin{lemma}[[Eq.~(17) and Lemma~4~\cite{lei18game}]\label{lem-recur}  For any $q\in (0,1)$ and $v>0$, the following hold:
\begin{align*}  &  \sum_{m=1}^{k+1 } q^{k+1 -m} m^{-v}    \leq  q^{k+1 } {e^{2v}q^{-1}-1 \over 1-q}+{2   (k+1)^{-v}\over q \ln(1/q) } , \quad \forall k \geq 0 \tag{i} ,\\
 \tag{ii}  & q^x   \leq c_{q,v}  {x^{-v}} \mbox{ for all $x >0$ where } c_{q,v} \triangleq   e^{- v  }{\left(\tfrac{v}{ \ln(1/q)}\right)^{ v  }}.
\end{align*}
\end{lemma}

\begin{proposition}[{\bf Polynomial rate of convergence}] \label{prp-rate-poly}  Let   Assumptions    \ref{ass-noise} and  \ref{ass-PL}    hold. Consider  Algorithm  \ref{alg-prox-gradient}, where   $p_i={1\over n} $, $ \alpha_i={2-\sqrt{3}\over  L_i}$ and  $N_i(k)= \Pi_{t=1}^v (\Gamma_i(k) +t)   $   for some
{positive integer $v\geq 1$.}
 Let  $C_v  \triangleq { (2-\sqrt{3})n^{v-1} \sigma^2 \over 2}  \sum_{i=1}^n  L_i^{-1}$,  $ \rho^* \triangleq  1- {(2-\sqrt{3})^2 \mu   \over n L_{\max}} $,  \blue{and $ C_f\triangleq
c_{\rho^*,v}  C_v{e^{2v}/\rho^*-1 \over 1-\rho^*}+c_{\rho^*,v}
\mathbb{E}\big[ F(x(0))-F^*\big] + {2 C_v \over \rho^*  \ln(1/\rho^*) }    $}. Then \begin{align} \label{rate-poly}
 \mathbb{E}\big[
F(x(k)) -F^*\big]   \leq   C_f k^{-v},  \quad \forall k \geq 1.
\end{align}
\end{proposition}
{\bf Proof.}
By setting  $\beta=\sqrt{3}-1$,  $ \alpha_i ={2-\sqrt{3} \over L_i}$, and $p_i={1\over n}$  in \eqref{x-bd3},   and using  $N_i(k)= \Pi_{t=1}^v (\Gamma_i(k) +t)   $  and \eqref{bd-batch}, we obtain the following recursion:
\begin{align} \label{bd-poly}
\mathbb{E}\big[ F(x(k+1)) -F^*\big] & \leq\rho^* \mathbb{E}\big[ F(x(k))-F^*\big] + (k+1)^{-v} { (2-\sqrt{3})\sigma^2 n^{v-1} \over 2}  \sum_{i=1}^n L_i^{-1} \notag
\\&   \leq (\rho^*)^{k+1}\mathbb{E}\big[ F(x(0))-F^*\big] + C_v \sum_{m=1}^{k+1} (\rho^*)^{k+1-m} m^{-v}.
\end{align}
This  together with Lemma \ref{lem-recur}(i)    {produces   the following inequality.  \begin{equation}\label{prp1-r2}
 \mathbb{E}\big[ F(x(k)) -F^*\big]   \leq    (\rho^*)^k  \left ( C_v{e^{2v}/\rho^*-1 \over  1-\rho^*}+\mathbb{E}\big[ F(x(0))-F^*\big] \right)+ {\tfrac{2 C_v k^{-v}}{\rho^*  \ln(1/\rho^* )}}  .
\end{equation}Since $\rho^* \in (0,1)$ and $v\geq 1$, by Lemma
\ref{lem-recur}(ii) we see that for any  $k\geq 1:$  $(\rho^*)^k \leq  c_{\rho^*,v}  k^{-v}   $  with $c_{\rho^*,v}   \triangleq e^{- v}{\left(\tfrac{v}{ \ln(1/\rho^*)}\right)^{ v  }} .$  This combined with  \eqref{prp1-r2} proves \eqref{rate-poly}. }
\hfill $\Box$
\subsection{Iteration and Oracle Complexity}

Next, we derive the iteration complexity (measured by the   number of proximal oracle) and oracle complexity (measured by the expected number of stochastic first-order oracle) bounds for obtaining an  $\epsilon-$optimal solution 
 such that
 $\mathbb{E} \left [F(x)\right] -F^* \leq \epsilon .$

\begin{theorem}\label{thm2} {\bf (Iteration Complexity)}
 Suppose   Assumptions  \ref{ass-noise} and  \ref{ass-PL}    hold.  Let  $\{ x(k)\}$ be generated by Algorithm \ref{alg-prox-gradient}, where   $ \alpha _i={2-\sqrt{3}\over  L_i},p_i={1\over n}$
and  $N_i(k)=\left \lceil (1-q_i)^{ -\Gamma_i(k) } \right \rceil$ for some $q_i\in (0, {(2-\sqrt{3})^2 \mu   \over   L_{i}})$.
Define  $q_{\min} \triangleq \min_{i\in \mathcal{N}} q_i$,   $q^*=1-{q_{\min}\over n}$,  and
$\eta_i \triangleq  { q_i\over n(1-q_i)}+1 $. Then   the  iteration   complexity and oracle complexity    required to obtain an   $\epsilon-$optimal solution are  respectively  $  \frac{
\ln( 1/\epsilonbar )}{\ln\left(1/ q^*\right)} $  and $ {1\over n} \sum_{i=1}^n{\eta_i\over  \ln (\eta_i)} (1/\epsilonbar)^{\frac{\ln(\eta_i)}{\ln\left(1/q^*  \right)}}  $ with $\bar{\epsilon}  $   defined by
\begin{align} \label{def-epsilon}
 &   \bar{\epsilon}  \triangleq   \epsilon \Big(    F(x(0))-F^*
+  { (1-\sqrt{3}/2)  \sigma^2\sum_{i=1}^n L_{\max}/L_i\over  {(2-\sqrt{3})^2 \mu   } -q_{\min} L_{\max}}    \Big)^{-1}.
\end{align}
  \end{theorem}
{\bf Proof.} Note from  $ q_{i}<  {(2-\sqrt{3})^2 \mu   \over   L_{i}}  $  for each $i\in \mathcal{N}$ that  $ q_{\min}<  {(2-\sqrt{3})^2 \mu   \over   L_{\max}} $.
Define $K_1(\epsilon)   \triangleq\left \lceil \frac{
\ln( 1/\epsilonbar)}{\ln\left(1/\rho^* \right)}\right \rceil$ with $ \bar{\epsilon} $   defined by \eqref{def-epsilon}.   Then    by \eqref{prp1-r1}, $ \alpha _i={2-\sqrt{3}\over  L_i},p_i={1\over n},$ and $ q^*=1-{q_{\min}\over n}$,    there holds: \begin{align*}
\mathbb{E}\big[ F(x_{ k }) -F^*\big] & \leq    \left(1-{q_{\min}\over n}  \right) ^{k} \left(     F(x(0))-F^*
+  { \sum_{i=1}^n \alpha_i \sigma^2/2 \over  {(2-\sqrt{3})^2 \mu   /   L_{\max}} -q_{\min}  }   \right)
\\&    = (q^*)^{k}  \epsilon   \bar{\epsilon}^{-1} \leq \epsilon \quad \forall k\geq  K_1(\epsilon) .
\end{align*}
Then the number of PO  to obtain an   $\epsilon-$optimal solution  is   $ K_1(\epsilon) $.
By  $p_i=1/n$ and \eqref{dis-ga},
 \begin{align}\label{rand-def-eta}
  \mathbb{E} \left[  N_i(k)  \right] &\leq   \mathbb{E} \left[(1-q_i)^{-\Gamma_{i ,k }}  \right]  +1
   = \sum_{m=0}^{k} \left({k \atop m}\right) \left (p_i(1-q_i)^{-1} \right)^m(1-p_i)^{k-m}  +1 \notag
 \\& =  (p_i(1-q_i)^{-1}+1-p_i)^{k}  +1=  \eta_i^{k}+1.
\end{align}
  Note   that for $\lambda > 1$, $
\sum_{k=0}^K \lambda^k  \leq \int_0^{K+1} \lambda^x dx \leq
\frac{\lambda^{K+1}} {\ln(\lambda)}. $
Since $ i_k $  is  independent of $N_i(k)$,
by using \eqref{def-epsilon} and \eqref{rand-def-eta}, the expected number of SFO calls required  to approximate an  $\epsilon-$optimal  solution  is   bounded as follows.
\begin{align*}&   \mathbb{E}\left[  \sum_{k=0}^{ K_1(\epsilon) -1} \sum_{i=1}^n N_i(k)I_{[i_k=i]} \right]
 =  \sum_{k=0}^{ K_1(\epsilon)-1 } \sum_{i=1}^n \mathbb{E}\left[ N_i(k) \right] \mathbb{E} \left[ I_{[i_k=i]}   \right]
\\ & =\sum_{k=0}^{ K_1(\epsilon)-1} \sum_{i=1}^n p_i \mathbb{E}\left[ N_i(k) \right]  \leq   \sum_{k=0}^{ K_1(\epsilon)-1}  {1\over n}\sum_{i=1}^n \left( \eta_i^{k}  +1  \right)={1\over n} \sum_{i=1}^n \sum_{k=0}^{ K_1(\epsilon)-1}  \eta_i^{k }   + K_1(\epsilon)
\\&  \leq  {1\over n} \sum_{i=1}^n { \eta_i   \over  \ln (\eta_i)} \eta_i^{  \frac{
\ln( 1/\epsilonbar)}{\ln\left(1/ q^*\right)} }+\left \lceil \frac{
\ln( 1/\epsilonbar)}{\ln\left(1/ q^*\right)} \right \rceil .\end{align*}
Note that for any $0<\epsilon, q<1$, we have the following relations:
$$  \eta^{  \frac{\ln( 1/\epsilon )}{\ln(1/q )} }
	  = \left(e^{\ln(\eta )}\right)^{  \frac{\ln( 1/\epsilon )}{\ln(1/q )} }
			  =  e^{\ln(1/\epsilon )^{\frac{\ln(\eta )}{\ln(1/q  )}}}
			  	=   {(1/\epsilon )^{\frac{\ln(\eta)}{\ln(1/q  )}}}  .$$
Thus,  the expected   number of  SFO  to obtain an    $\epsilon-$optimal  solution  is bounded by
 \begin{align}\label{bd-m}
 {1\over n} \sum_{i=1}^n { \eta_i   \over  \ln (\eta_i)}{(1/\epsilonbar)^{\frac{\ln(\eta_i)}{\ln\left(1/q^* \right)}}}  +\left \lceil \frac{\ln( 1/\epsilonbar)}{\ln\left(1/q^* \right)} \right \rceil ,
\end{align}
giving us the required oracle complexity.
\hfill $\Box$

The following corollary   emerges from Theorem \ref{thm2} with  $q_i$ taking a specific form.
\red{It shows  how  the number of blocks $n$, the Lipschitz constants, and the parameter  $\mu$ in
 Assumption  \ref{ass-PL}  collectively   influence the constants as well as the order  in the iteration and oracle complexity. }
\begin{corollary}
Let $N_i(k)=\left \lceil (1-q_i)^{ -\Gamma_i(k) } \right \rceil$ with $q_i = \frac{\theta_i \mu}{L_i}$ for some  $\theta_i \in (0,(2-\sqrt{3})^2)$ in Algorithm \ref{alg-prox-gradient}, and  the other conditions of Theorem \ref{thm2}  still hold. Define $\theta_{\min}\triangleq \min_{i } \theta_i$ and $\theta_{\max}\triangleq \max_{i } \theta_i$.
Then the  iteration and oracle complexity required  to obtain an
$\epsilon-$optimal solution are   $\mathcal{O}
\left(\frac{nL_{\max}}{\mu} \ln(1/\epsilon)\right)$
and $\mathcal{O}\Big({n L_{\rm ave } \over \mu} (1/\epsilon)^{  \big(1+{1\over n\kappa_{\min}-1}  \big) \big({L_{\max}  \theta_{\max} \over L_{\min}  \theta_{\min} } \big) }\Big)$, respectively.
\end{corollary}
{\bf Proof.}  We begin by deriving a bound on $K_1(\epsilon)$:
\begin{align*}
\frac{\ln(1/\bar{\epsilon})}{\ln(1/q^* )} =
\frac{\ln(1/\bar{\epsilon})}{-\ln(q^*)} =
\frac{\ln(1/\bar{\epsilon})}{-\ln(1-q_{\min}/n)} \leq
\frac{\ln(1/\bar{\epsilon})}{q_{\min}/n} = \frac{nL_{\max}}{ \theta_{\min}  \mu} \ln(1/\bar{\epsilon}),
\end{align*}
where $-\ln(1-\tfrac{q_{\min}}{n}) \geq \tfrac{q_{\min}}{n}$ and  $q_{\min}  \geq  \tfrac{ \theta_{\min} \mu}{L_{\max}}$,  implying that $ K_1(\epsilon)=\mathcal{O}\left(\frac{nL_{\max}\ln(\tfrac{1}{\epsilon})}{\mu}\right).$
 Next, we analyze the two terms necessary for bounding the oracle complexity.
\begin{align*}
&\frac{\ln(\eta_i)}{\ln(1/q^* )}   \leq
\frac{nL_{\max}}{\theta_{\min} \mu} \ln(\eta_i) =
\frac{nL_{\max}}{ \theta_{\min}  \mu} \ln(1+q_i/(n-q_i)) \leq  \frac{nL_{\max}}{\theta_{\min} \mu} \frac{q_i}{ n-q_i }   \\
& \leq
\frac{nL_{\max}}{\theta_{\min} \mu} \frac{\mu \theta_{\max}}{L_i (n-\theta_{i} \mu /L_i)}
 \leq   \frac{nL_{\max} \theta_{\max}}{  (n-  \mu /L_{\min}) L_{\min} \theta_{\min}}=  \frac{nL_{\max} \theta_{\max}}{  (n-  1 /\kappa_{\min}) L_{\min} \theta_{\min}},
\end{align*}
 where the second inequality holds by   $\ln(1+x) \leq x , \forall x \in [0,1),$  the fourth inequality follows from $\theta_i \leq 1$, and  the last equality follows from
 $\kappa_{\min}=L_{\min} /\mu.$
 In addition, we derive a bound on $\eta_i/\ln(\eta_i)$  where $\eta_i = 1+\tfrac{q_i}{n(1-q_i)}$.
 Since $\ln(1+x) \geq x/(x+1)$ for any $x\geq 0$,  we then  have
\begin{align*}
	\frac{\eta_i}{\ln(\eta_i)} & = \frac{ 1+q_i/(n(1-q_i)) }{\ln(1+q_i/(n(1-q_i)))}
	 = \frac{(1+x^0)}{\ln(1+x^0)} \leq \frac{(x^0+1)^2}{x^0}  \\
&  =  x^0+2+{1\over x^0}= 2 + \frac{q_i}{n(1-q_i)} + {n(1-q_i) \over q_i}\leq 2 + \frac{1}{n}+ \frac{n L_i}{\theta_i \mu}   ,
\end{align*}
 where $x^0=q_i/(n(1-q_i))$ and $q_i/(1-q_i) \leq 1$ (since $q_i \leq (2-\sqrt{3})^2$).
 We prove the result by deriving a bound on the oracle complexity:
\begin{align*}
& \quad {1\over n} \sum_{i=1}^n { \eta_i   \over  \ln (\eta_i)}{(1/\epsilonbar)^{\frac{\ln(\eta_i)}{\ln\left(1/\rho^*  \right)}}}  +\left \lceil \frac{\ln( 1/\epsilonbar)}{\ln\left(1/\rho^*  \right)} \right \rceil
\\& {{{\scriptstyle   \eqref{def-epsilon} }\atop=}   \mathcal{O}\Bigg( {\sum_{i=1}^n L_i \over \mu} \left({\mathbb{E}[F(x(0))]-F^*\over \epsilon}+ {n \sigma^2\over \epsilon} \right)^{ \frac{nL_{\max} \theta_{\max}}{  (n-1 /\kappa_{\min})L_{\min} \theta_{\min}}} \Bigg)  }.  \qquad\qquad\qquad \Box
\end{align*}
It can be seen \red{from the above corollary} that if $\theta_{\max} = \theta_{\min}$ (by choosing $\theta_i = \theta$ for all $i$) and $L_{\max} = L_{\min}=L$,   the oracle complexity  becomes  $\mathcal{O}({n\kappa }(1/\epsilon)^{1+ { 1 \over  n\kappa -1}})$
with $\kappa\triangleq {L \over \mu}$, which tends to the  optimal oracle complexity of $\mathcal{O}\left ({n \kappa \over \epsilon} \right)$ for large $n.$
From Theorem~\ref{thm2}, we may obtain the {\em optimal} oracle complexity for  $n = 1$. \red{This is because
$\ln(\eta_i)/\ln(1/\rho^*  ) = \ln(1/(1-q))/\ln(1/(1-q)) = 1 $ when $n=1.$}

\begin{corollary}\label{cor-single}
 Suppose $n=1 $, Assumptions  \ref{ass-noise} and  \ref{ass-PL}    hold.  Consider  Algorithm \ref{alg-prox-gradient}, where $ \alpha={2-\sqrt{3}\over  L}$ and  $N(k)=\left \lfloor (1- q)^{ -k} \right \rfloor$ with $ q\in \left (0,   {(2-\sqrt{3})^2  \mu  \over  L }\right) $.
Then the  iteration and oracle complexity required to obtain an   $\epsilon-$optimal solution are   $\mathcal{O}({L\over \mu}\ln(1/\epsilon))$  and  $\mathcal{O}\left({L\over \mu \epsilon} \right)$ ,  respectively.
\end{corollary} 

In Theorems \ref{prp1} and   \ref{thm2}, we establish the rate as well as the   iteration and  oracle complexity
bounds of   Algorithm \ref{alg-prox-gradient} when each block   is randomly picked with equal probability.
When blocks  are chosen by a non-uniform distribution, we state a
result in the smooth regime but omit the proof since  it
is similar to   Theorems \ref{prp1} and \ref{thm2}.

\begin{corollary}\label{cor} Suppose $r(x)\equiv 0$ and  $\bar{f}$
satisfies $ \|\nabla \bar{f}(x)\|^2 \geq 2\mu
\big(\bar{f}(x)-F^*\big)$ with $\mu>0.$   Consider Algorithm \ref{alg-prox-gradient}, where   $
\alpha_i={2-\sqrt{3}\over  L_i}$, $p_i={L_i\over \ \sum_{j=1}^n L_j}
$, and $N_i(k)=\left \lceil (1-q_i)^{ -\Gamma_i(k) } \right \rceil$
with $q_i = \frac{\theta_i \mu}{L_i}$ for some $\theta_i \in
(0,(2-\sqrt{3})^2)$.   Then the  iteration and oracle complexity bounds for obtaining an   $\epsilon-$optimal solution are
$\mathcal{O}\left({n L_{\rm ave}\over \mu}\ln \left({1\over
\epsilon} \right) \right)$  and  $\mathcal{O}\Big( {n L_{\rm ave }
\over \mu }(1/\epsilon)^{\big(1+{1\over n\kappa_{\min}-1}  \big)
\frac{ L_{\rm ave} \theta_{\max}}{ L_{\min} \theta_{\min}}}\Big)$,
respectively.  \end{corollary}

\red{This result shows that the non-uniform block selection rule can   improve the constant in the iteration complexity when compared with the uniform block selection.}
 Finally,  we investigate the iteration and oracle complexity of  Algorithm \ref{alg-prox-gradient} when   the  number of the sampled gradients increases at a  slower polynomial  rate.
\begin{proposition} \label{prp-complexity } Let Assumptions  \ref{ass-noise} and  \ref{ass-PL}    hold. Consider Algorithm \ref{alg-prox-gradient}, where  $p_i={1\over n}$,  $ \alpha_i={2-\sqrt{3}\over  L_i}$,  and  $N_i(k)= \Pi_{t=1}^v (\Gamma_i(k) +t)   $ for some     integer $v\geq 1$.
Define $ \rho^* \triangleq  1- {(2-\sqrt{3})^2 \mu   \over n L_{\max}} $ and $C_v  \triangleq { (2-\sqrt{3})n^{v-1} \sigma^2 \over 2}  \sum_{i=1}^n  L_i^{-1}$.
Then the  iteration and oracle complexity bounds required to obtain an  $\epsilon-$optimal solution  are $\mathcal{O} (v (1/\epsilon)^{1/v})$ and   \blue{$\mathcal{O} \left(e^v v^{2v+1} (1/\epsilon)^{1+1/v}  \right)$}, respectively.
\end{proposition}
{\bf Proof.}  From \eqref{rate-poly} it follows that for any    $ k\geq  K(\epsilon)\triangleq  \left( {
C_f \over \epsilon }\right)^{1/v}$,  $ \mathbb{E}\big[ F(x(k))
-F^*\big]    \leq \epsilon.$  \blue{From the definition of $C_f$ in Proposition   \ref{prp-rate-poly} it follows that}   $C_f=\mathcal{O}(e^v v^v)$,   hence the iteration complexity is  $\mathcal{O} ( v (1/\epsilon)^{1/v})$ .
Note that $(m+v)^v =\sum_{t=0}^v  \left({v \atop t }\right)   m^t v^{ v-t}.$ Then by $N_i(k)=\Pi_{t=1}^v (\Gamma_i(k) +t)  \leq  (\Gamma_i(k) +v)^v$  and  using \eqref{bd-Nk},  we obtain that for each   $    i \in \cal{N}:$
\vspace{-0.2in}
 \begin{equation}\label{expect-nk}
\begin{split}
 & \mathbb{E} \left[  N_i(k)   \right]     \leq \sum_{m=0}^{k} (m+v)^v \mathbb{P}(\Gamma_i(k)=m)
  \\& = \sum_{t=0}^v  \left({v \atop t }\right)   v^{v-t} \sum_{m=0}^{k}   m^t  {k!    \over m! (k-m)!  }  p_i^{m } (1-p_i)^{k-m} = \sum_{t=0}^v  \left({v \atop t }\right)   v^{v-t}  \mathbb{E}[\Gamma_i(k)^t] .
\end{split}
\end{equation}
By \cite[p.154]{papoulis1984probability} we know  that the $t$-th moment of the  binomial distribution $\Gamma_i(k) $
equals the $t$-th  derivative of $M_i(y)$  at $y=0 $, where  $M_i(y)=(p_i e^y+1-p_i)^k $. Thus, we can   show that
 $  \mathbb{E}[\Gamma_i(k)^t]=\mathcal{O}(k^t) $, hence by \eqref{expect-nk} and $  \left({v \atop t }\right)   v^{v-t}=\mathcal{O}(v^v)$  we obtain  $ \mathbb{E} \left[  N_i(k)   \right]  =\blue{\mathcal{O}(v^{v+1} k^v)}.$  Note
 \blue{by    $C_f=\mathcal{O}(e^v v^v)$} that
$$  \sum_{k=0}^{ K(\epsilon) -1 }  k^v \leq \int_{1}^{ K(\epsilon)  } t^v dt= {t^{v+1} \over v+1}\Big |_{1}^{ K(\epsilon) }\leq  (v+1)^{-1} \left({ C_f \over \epsilon }\right)^{1+{1\over v}}\blue{=\mathcal{O} (v^ve^v(1/\epsilon)^{1+1/v})) .}$$
 \blue{This together with $ \mathbb{E} \left[  N_i(k)   \right]  =\mathcal{O}(v^{v+1} k^v)$ implies that}  the expectation of the total number of sampled gradients   required to obtain   an   $\epsilon-$optimal solution   is   bounded by
 $$\mathbb{E}\Big[  \sum_{k=0}^{ K (\epsilon) -1} \sum_{i=1}^n N_i(k)I_{[i_k=i]} \Big]
  =\sum_{k=0}^{ K (\epsilon)-1} \sum_{i=1}^n p_i \mathbb{E}\left[ N_i(k) \right]  =\blue{\mathcal{O} (v^{2v+1}e^v (1/\epsilon)^{1+1/v})}.    \qquad \qquad \Box$$

\blue{
 \begin{remark}\label{rem-poly}  Proposition \ref{prp-complexity }  implies that as the polynomial degree $v$ is increased, the constant in   the iteration and oracle complexity  respectively grows at a linear and an exponential rate, respectively. Therefore, choosing an appropriate $v$ requires assessing both available computational resources and the impact of   generating   sample-average gradients  with large sample-sizes.
\end{remark}}

\section{Numerical Experiments}\label{Sec:exa}

In this section, we will  examine  empirical  algorithm performance   on the  sparse and nonlinear least squares problems to demonstrate the behavior of  Algorithm \ref{alg-prox-gradient}.
 Throughout this section, the  empirical mean error is based on averaging across  50 trajectories.

\subsection{ Sparse  Least Squares}
We  apply Algorithm \ref{alg-prox-gradient} to the following sparse  least squares   problem:
\begin{align}\label{exp1}
\min_{ x\in \mathbb{R}^{  d}}{1\over 2N}\sum_{i=1}^N ( a_i^T x-b_i)^2+\lambda \| x\|_1, \quad \lambda >0.  \tag{LASSO}
\end{align}
 We  first generate a sparse vector $x^*$  where  only $10\%$  components of the vector is nonzero with nonzero ones independently   generated from the standard normal distribution. We  then generate $N$ samples
$(a_i,b_i)$, where  components of  $a_i\in \mathbb{R}^d$ are generated from standard normal distribution while
 $b_i= a_i^Tx^*+ \hat{\epsilon}$, where  $\hat{\epsilon}$ is normally distributed  with zero mean    and standard deviation   $0.01$. We partition $x \in \mathbb{R}^d$ into $n=10$ blocks and set $\lambda=0.1$.

{\bf Sensitivity to sample-size policies:}  We now  implement  Algorithm
\ref{alg-prox-gradient} with  $\alpha=0.01, ~p_i={1\over n}$,  the geometric
batch-size $N_i(k)=\left \lceil q^{ -\Gamma_i(k) } \right \rceil$, and
investigate  how the  parameters    $q,N,d$ influence  the algorithm
performance.  We ran Algorithm \ref{alg-prox-gradient}  for  50 epochs where
each epoch  implies the usage of all samples.   The results are displayed  in Table
\ref{TAB2}  for  the empirical relative  error $\mathbb{E}[F(x )] -F^*
\over F^*$, the number of proximal evaluations, and CPU times.  The
results  suggest   that for given a fixed simulation budget,   slower   geometric
rates of growth of  batch-sizes  lead to better empirical error while
requiring more CPU time since more proximal evaluations are needed.  In
addition,  it is noticed that the  running  time   increases approximately  linearly  with   $N$
and  $d $.

\begin{table}[htbp]
 \small
    \begin{minipage}[c]{0.48 \textwidth}
\centering  (a) $d=400$ \\ { \begin{tabular}{|c|c|c|c|c| }
 \hline $N$ & $q$&   emp.err  &
 prox.eval  & CPU(s)   \\ \hline  \hline
{\multirow{3}{*}{$ 1000$}}  &  $ 0.85$ & 2.46e-02 & 86&  {4.22}   \\   \cline{2-5}
 &  $  0.9 $     &   1.71e-02 & 105 &4.7  \\ \cline{2-5}
 &  $  \bf{0.95}$    &\bf{5.00e-03} & 164  &  {6.28}      \\ \hline
{\multirow{3}{*}{$ 2000$}}  &  $0.85$ &  3.71e-02  & 90 &  {7.71}   \\   \cline{2-5}
 &  $  0.9 $     &   2.49e-02   & 112&8.82 \\ \cline{2-5}
 &  $  \bf{0.95}$    &\bf{6.10e-03 } & 178  &  {11.48}   \\ \hline
{\multirow{3}{*}{$ 4000$}}  &  $ 0.85$ &  1.27e-02  & 94 & {16.15}  \\   \cline{2-5}
 &  $  0.9 $     &  7.60e-03  & 119&18.27 \\ \cline{2-5}
 &  $  \bf{0.95}$    &\bf{1.90e-03} & 192  &  {24.3}   \\ \hline
 \end{tabular}
}
\end{minipage}
 \begin{minipage}[c]{0.48 \textwidth}
\centering
(b) $d=800$ \\{ \begin{tabular}{|c|c|c|c|c|}
 \hline  $N$ &$q$ &   emp.err  &
 prox.eval  & CPU(s)   \\ \hline  \hline
{\multirow{3}{*}{$ 1000$}}  &  $ 0.85$ &  2.8e-02  & 86&  {8.57}   \\   \cline{2-5}
 &  $  0.9 $     &   1.93e-02 & 105 &9.69  \\ \cline{2-5}
 &  $  \bf{0.95}$    &\bf{7.10e-03} & 164&  {12.33}      \\ \hline
{\multirow{3}{*}{$ 2000$}}  &  $ 0.85$ &  1.62e-02  & 90 &  {17.28}   \\   \cline{2-5}
 &  $  0.9 $     &   1.10e-02  & 112&17.7 \\ \cline{2-5}
 &  $  \bf{0.95}$    &\bf{3.70e-03 } & 178  & {25.9}   \\ \hline
{\multirow{3}{*}{$ 4000$}}  &  $ 0.85$ &  1.62e-02 & 94 &  {34.2}  \\   \cline{2-5}
 &  $  0.9 $     &  1.00e-02   & 119&38.6 \\ \cline{2-5}
 &  $  \bf{0.95}$    &\bf{2.6e-03} & 192  &  {50.37}   \\ \hline
 \end{tabular}
}
\end{minipage} \vskip 2mm
      \caption{Comparison of the different selections of  batch-sizes \label{TAB2}}
\end{table}

{\bf Comparison  with  BSG  \cite{xu2015block}:} Let $N=2000 $ and $d=200$ in \eqref{exp1}.
We compare  Algorithm \ref{alg-prox-gradient}  with  BSG \cite{xu2015block}  by
running  both schemes for 50 epochs.  We show the results in  Table~\ref{tab3}
and  plot trajectories  in \red{Figures \ref{fig_BSG_iter} and \ref{fig_BSG_oracle}},  where   BSG-$m$ denotes the
minibatch BSG   that  utilizes $m$ samples  at each iteration while in
 Algorithm  \ref{alg-prox-gradient}, we set $N_i(k)=\left \lceil q^{ -\Gamma_{i
,k } } \right \rceil$.  The  empirical   rate of convergence   in terms of proximal evaluations  shown in \red{Figure \ref{fig_BSG_iter}}  implicitly supports  the
iteration complexity statements.  We  observe the following: (i) at first,
minibatch BSG  {displays a faster decay in objective} than Algorithm \ref{alg-prox-gradient} since the
batch-size in our scheme is relatively  small    at the   outset;
(ii)   Algorithm  \ref{alg-prox-gradient}   proceeds to catch  up and  outperform the  minibatch BSG
since the variance of the sampled gradient decreases with increasing
batch-size; (iii) Both   minibatch BSG with larger batch-sizes and  Algorithm  \ref{alg-prox-gradient} with
faster increasing batch-size  display   faster empirical rates with   fewer proximal
evaluations.  The  empirical  algorithm performance  in terms
of epochs     displayed  in \red{Figure \ref{fig_BSG_oracle}} demonstrates the results of oracle complexity. By comparing  the number of samples  given the  fixed  relative error,  Algorithm  \ref{alg-prox-gradient}
with  $N_i(k)=\left \lceil 0.98^{ -\Gamma_i(k) } \right \rceil$ has the
best performance, which can also be concluded   from Table~\ref{tab3}.

\begin{figure}[!htb]
    \centering
    \begin{minipage}{.48\textwidth}
        \centering
      \includegraphics[width=2.5in]{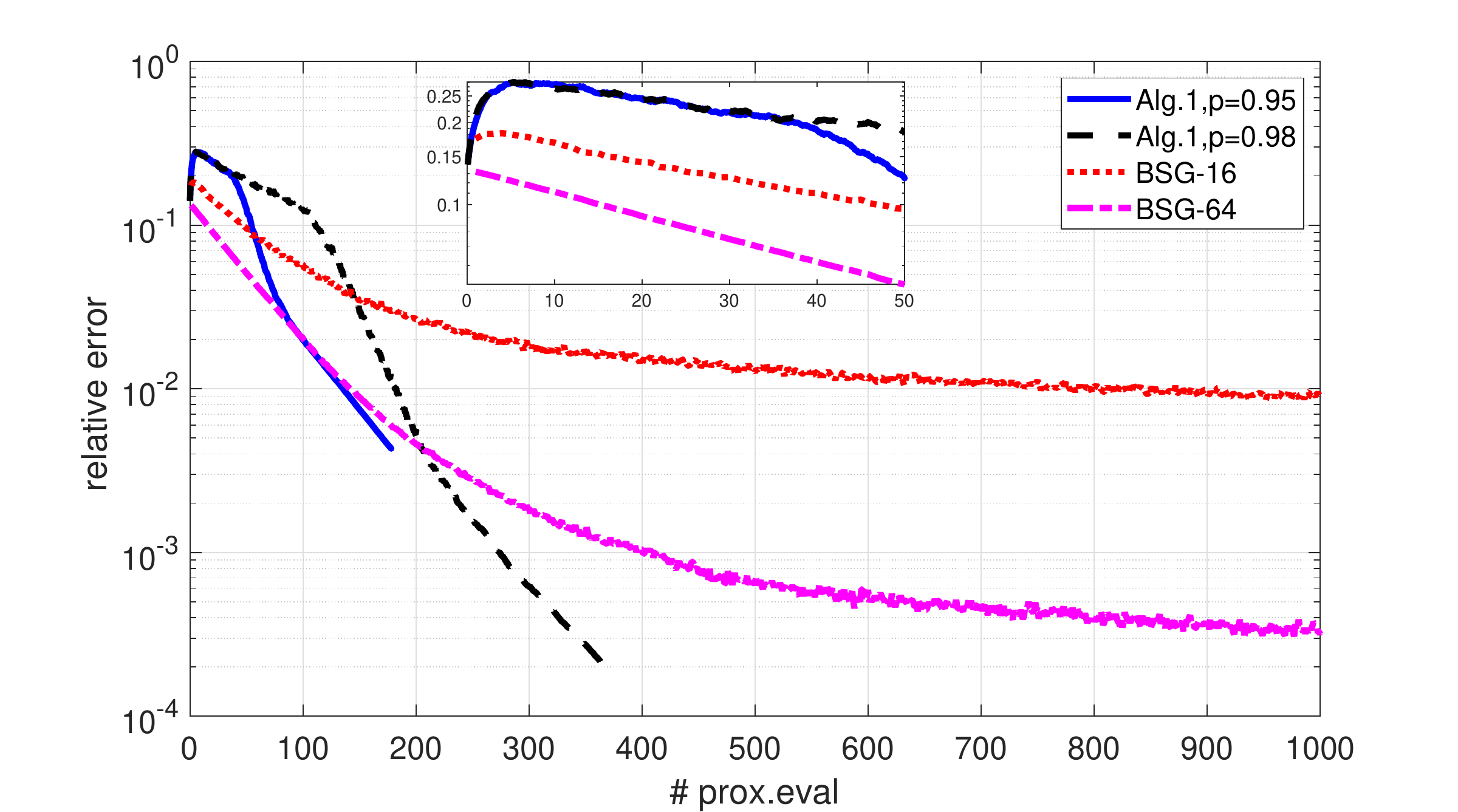}
\caption{ Iteration Complexity   of  Algorithm\ref{alg-prox-gradient} and  BSG    } \label{fig_BSG_iter}
    \end{minipage}%
    \begin{minipage}{0.48\textwidth}
        \centering
    \includegraphics[width=2.5in]{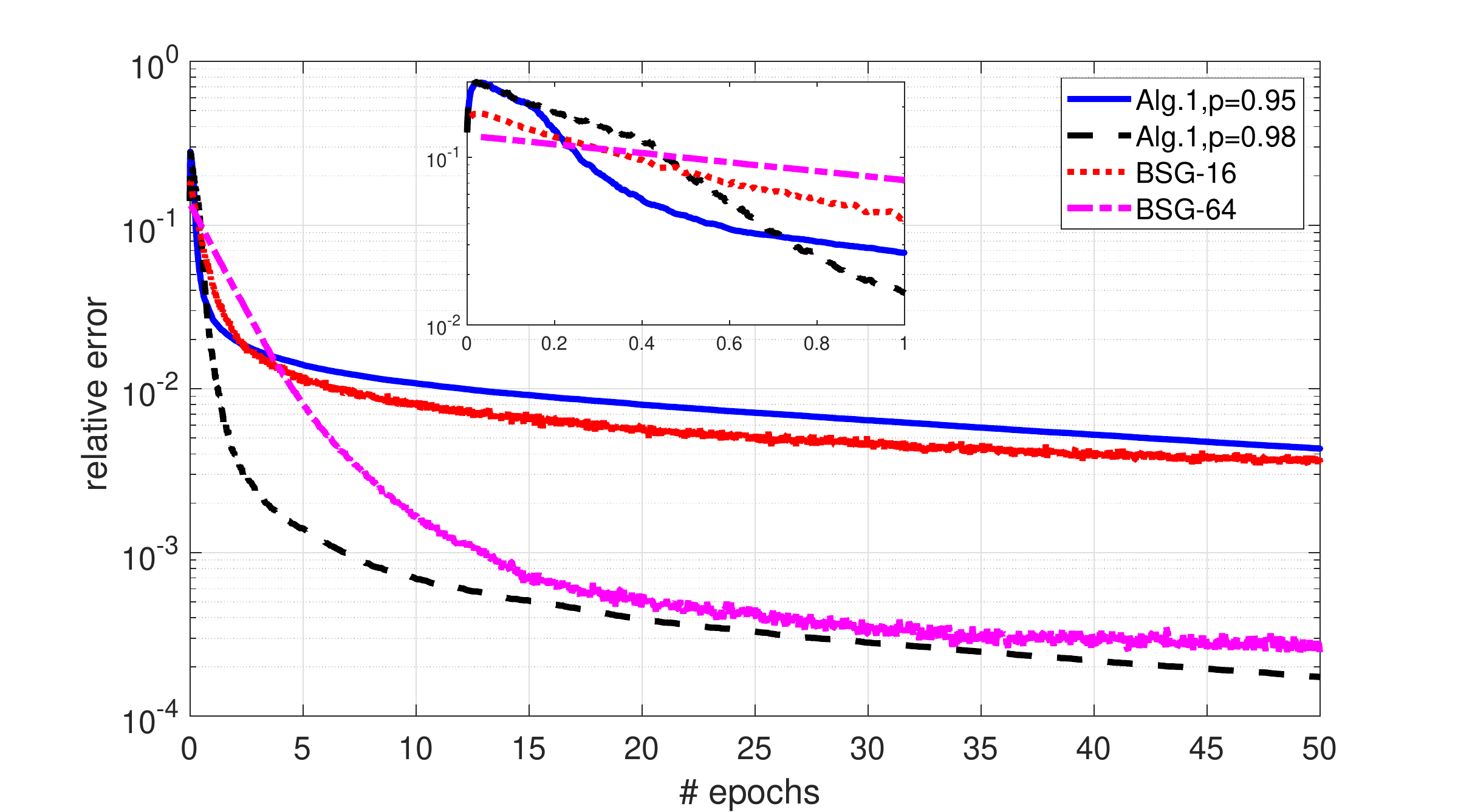}
\caption{Oracle Complexity   of  Algorithm\ref{alg-prox-gradient} and  BSG } \label{fig_BSG_oracle}
    \end{minipage}
\end{figure}

\begin{table}[htbp]
 \small
\newcommand{\tabincell}[2]{\begin{tabular}{@{}#1@{}}#2\end{tabular}}
\centering
 \begin{tabular}{|c|c|c|c|c| }
 \hline    &     \tabincell{c} { Algorithm  \ref{alg-prox-gradient}, p=0.95}    &  \tabincell{c} { Algorithm   \ref{alg-prox-gradient}, p=0.98} & BSG-16&  BSG-64
 \\ \hline  \hline
 emp.err& 4.30e-3 & \bf{1.73e-4} &2.60e-3 &2.75e-4\\ \hline
  prox.eval &   178 & 375&6251  &1563 \\ \hline
 CPU(s)   & 5.94 & \bf{10.25} &118.65&33.36 \\ \hline
 \end{tabular}
  \vskip 2mm \caption{Comparison of  Algorithm \ref{alg-prox-gradient}  and  BSG \label{tab3}}
\end{table}

 {\bf  Influence of  block-specific steplengths}:   In this experiment,   we set
$N=1000, d=200$, and  let  the   entries of  $a_i\in \mathbb{R}^d$ corresponding
to different blocks  be generated from normal distributions with zero mean but
with \uvs{differing}  variances. Such data generation implies that  the block-wise Lipschitz constants of  ${1\over 2N}\sum_{i=1}^N ( a_i^T x-b_i)^2$ \uvs{can vary}. We implement  Algorithm \ref{alg-prox-gradient}    with  the non-uniform block selection as per a distribution $p_i= {  L_i \over \sum_{i=1}^nL_i}$ in two settings: (i)  the same steplength $\alpha_i\equiv \alpha= {1.28 \over L}$
depending on the   Lipschtiz constant of $\nabla \bar{f}(x)$,  and (ii) the
block-specific steplength   $\alpha_i={1\over L_i}$  depending on   the
block-wise Lipschitz constant $L_i.$   Such a  selection ensures that the
steplengths in   two settings are approximately  the same   when the block-wise
Lipschitz constants are identical. For a particular set of realizations with
the Lipschtiz constant   satisfying $L_{\max}/L_{ave}=2.5$,  the empirical
iteration and oracle complexity of   Algorithm \ref{alg-prox-gradient} in the
two settings are shown in  Figures
\ref{fig_lipschitz} and \ref{fig_lipschitz0}, respectively.  These findings  reinforce the point that
block-specific steplengths, reliant on block-wise Lipschitz constant $L_i$,
display better empirical behavior since less proximal evaluations
(see Figure \ref{fig_lipschitz}) and less sampled gradients (see Figure
\ref{fig_lipschitz0}) are required  for obtaining a solution with  similar  accuracy.   In
addition,   we generate  four  sets of data, for which    the global
Lipschitz constant $L$  of the problem \eqref{exp1}  is the same while  the
ratio  $L_{\max}/L_{ave}$ \uvs{is} different. We  then  run  Algorithm
\ref{alg-prox-gradient}  with the  identical  and block-specific steplengths
on the four  generated datasets  up to  100 epochs and compare the  empirical
errors. The results are shown   in Table \ref{tab5}, where  $x^I(K)$ and
$x^B(K)$ denote  the estimates  generated by   Algorithm
\ref{alg-prox-gradient}   with   identical  and block-specific steplengths,
respectively.  Since the ratio    ${ \mathbb{E}[F(x^I(K) )] -F^* \over
\mathbb{E}[F(x^B(K) )] -F^* }$  is greater than one,  we  \uvs{may}  conclude that
the block-specific steplengths  might lead to much better algorithm performance
compared  with   identical   steplength.  We observe that empirical error can be $50$ times poorer when $L_{\max}/L_{ave} = 1.47$.

\begin{figure}[!htb]
    \centering
    \begin{minipage}{0.1\textwidth}
        \end{minipage}
    \begin{minipage}{0.43\textwidth}
        \centering
    \includegraphics[width=2.5in]{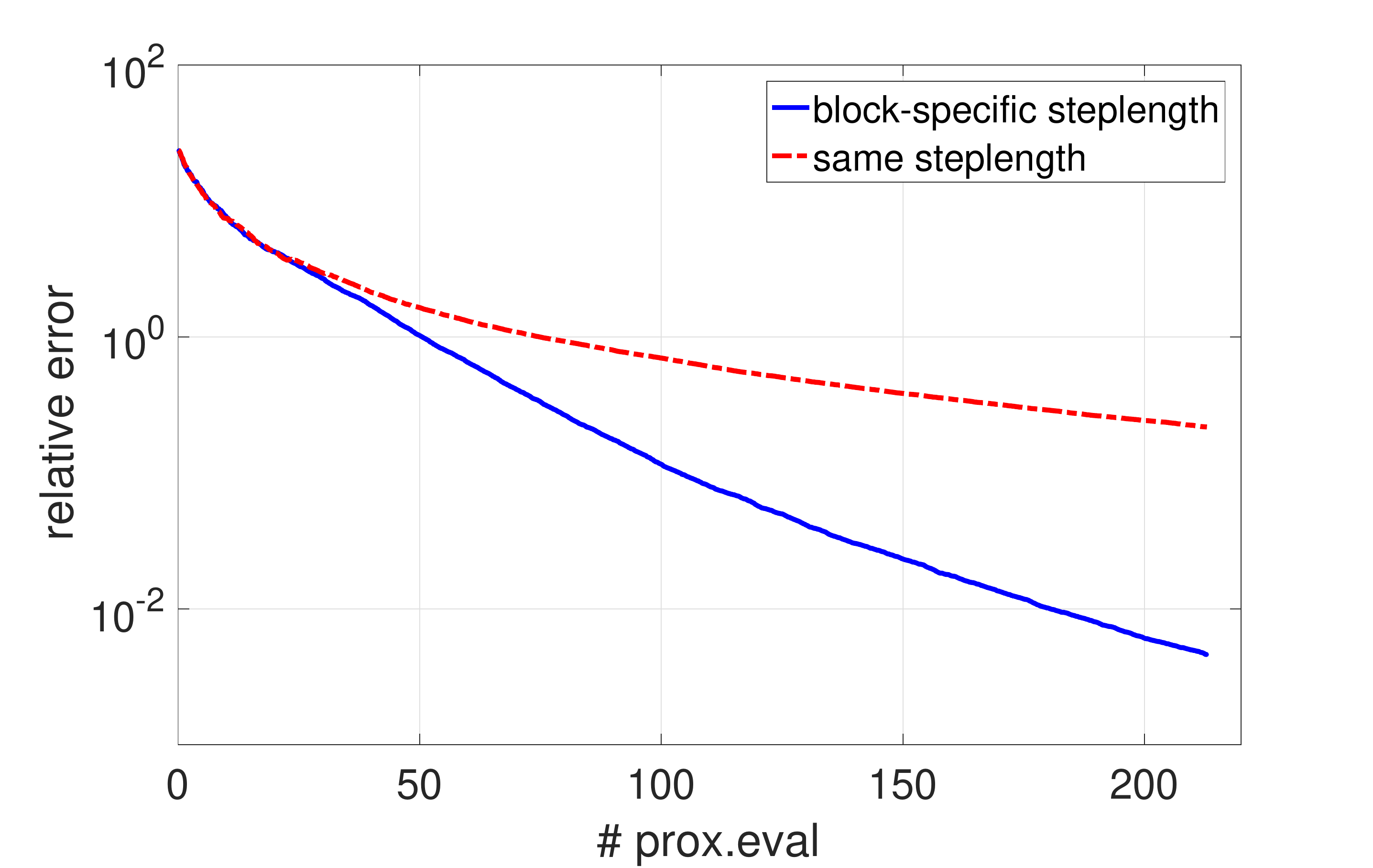}
\caption{Iteration Complexity of Algorithm \ref{alg-prox-gradient}  with  identical and block-specific steplengths   } \label{fig_lipschitz}
    \end{minipage}
    \begin{minipage}{.43\textwidth}
        \centering
      \includegraphics[width=2.5in]{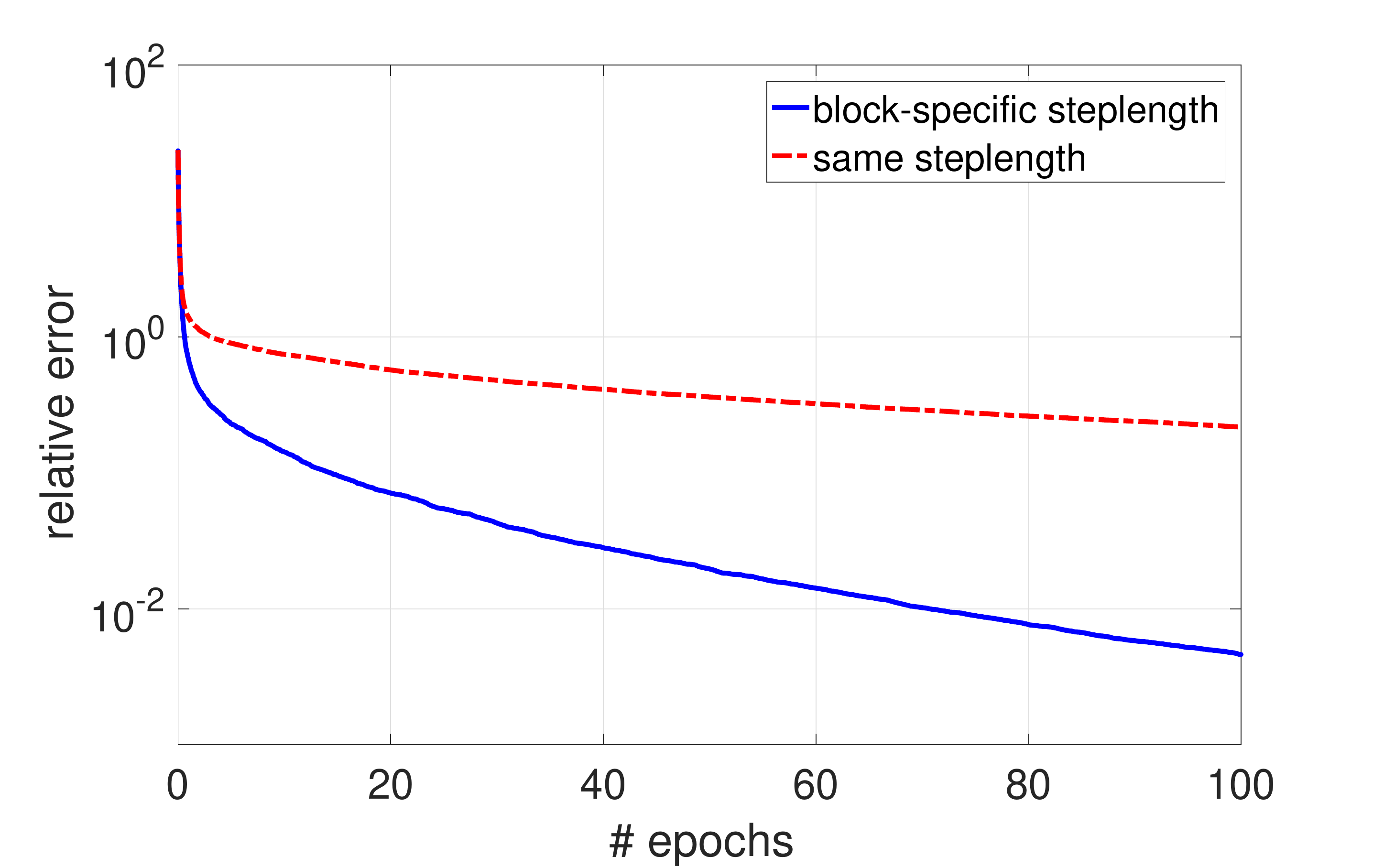}
\caption{Oracle Complexity of Algorithm \ref{alg-prox-gradient}  with identical and block-specific steplengths   }
\label{fig_lipschitz0}
    \end{minipage}%
\end{figure}

\begin{table}[htbp]
\centering
 \begin{tabular}{|c|c|c|c|c| }
 \hline    &     ${L_{\max} \over L_{ave}}=1.15$     &${L_{\max} \over L_{ave}}=1.27$
&${L_{\max} \over L_{ave}}=1.34$ &${L_{\max} \over L_{ave}}=1.47$
 \\  \hline
  ${ \mathbb{E}[F(x^I(K) )] -F^* \over \mathbb{E}[F(x^B(K) )] -F^* }$ & 15.3 & 27.5  & 31.9& 52.4 \\ \hline
 \end{tabular}
  \caption{Comparison of  Algorithm \ref{alg-prox-gradient}   with   identical  and block-specific steplengths  \label{tab5}}
\end{table}



\blue{ {\bf Uniform vs non-uniform block selection mechanism. } We  generate a particular set of realizations with
 Lipschitz  constants  of  ${1\over 2N}\sum_{i=1}^N ( a_i^T x-b_i)^2$
satisfying $L_{\max}/L_{ave}=1.35$. We implement  Algorithm
\ref{alg-prox-gradient}    with  the non-uniform  and the uniform block
selection probability, while  keeping  the other algorithm  parameters the
same. The empirical iteration and oracle complexities   are displayed in
Fig. \ref{fig_block}, showing that non-uniform  block selection
leads to   better performance than the uniform case and reduces  the
empirical error by approximately 50\%.  }

\begin{figure}[!htb]
    \centering
    \begin{minipage}{.46\textwidth}
        \centering
      \includegraphics[width=2.8in]{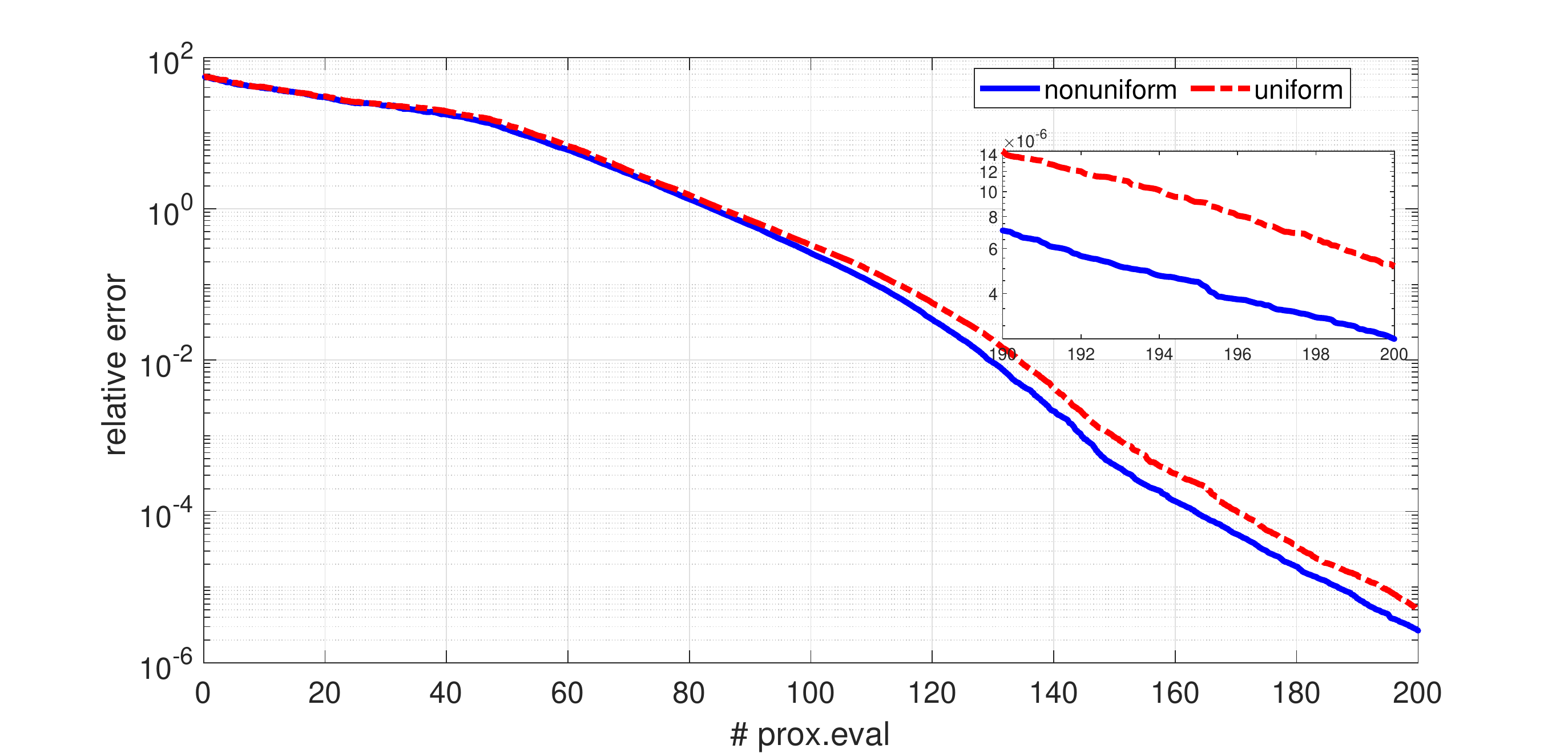}
    \end{minipage}%
    \begin{minipage}{0.46\textwidth}
        \centering
    \includegraphics[width=2.8in]{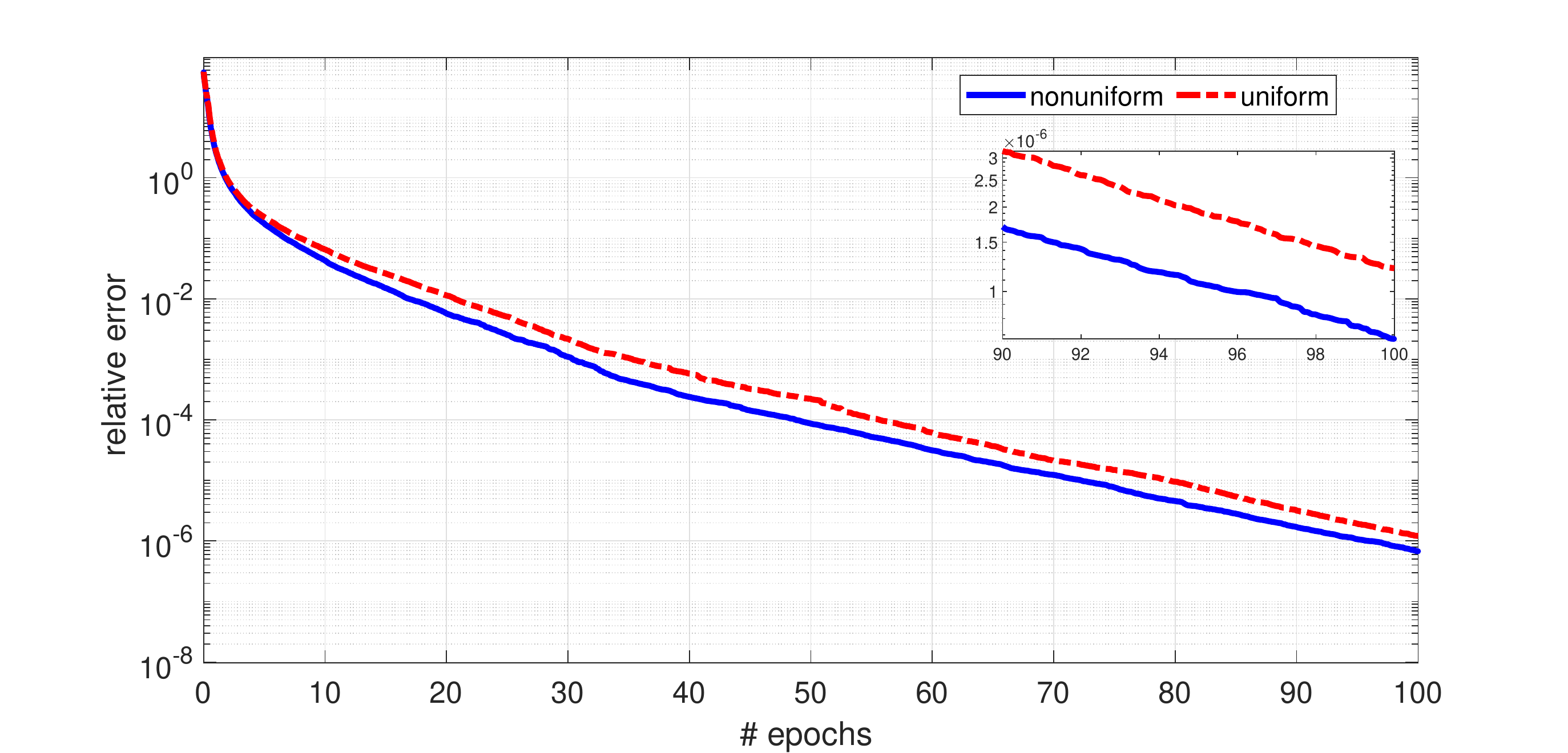}
    \end{minipage}
    \caption{   Algorithm~\ref{alg-prox-gradient} with Uniform and Non-uniform Block Selection  for \eqref{exp1}  with $L_{\max}/L_{ave}=1.35$} \label{fig_block}
\end{figure}

\blue{{\bf Influence of delays.}  We now compare the empirical  performance of  Algorithm \ref{alg-prox-gradient}    with   and without delays, where in the delayed case, we set the delay to be uniformly bounded by some  positive  integer.  The   iteration and oracle complexity  bounds are displayed in Fig. \ref{fig_delay}, which shows that
a larger delay leads to worse performance than the case with smaller delay or without delay.    }

\begin{figure}[!htb]
    \centering
    \begin{minipage}{.46\textwidth}
        \centering
      \includegraphics[width=2.8in]{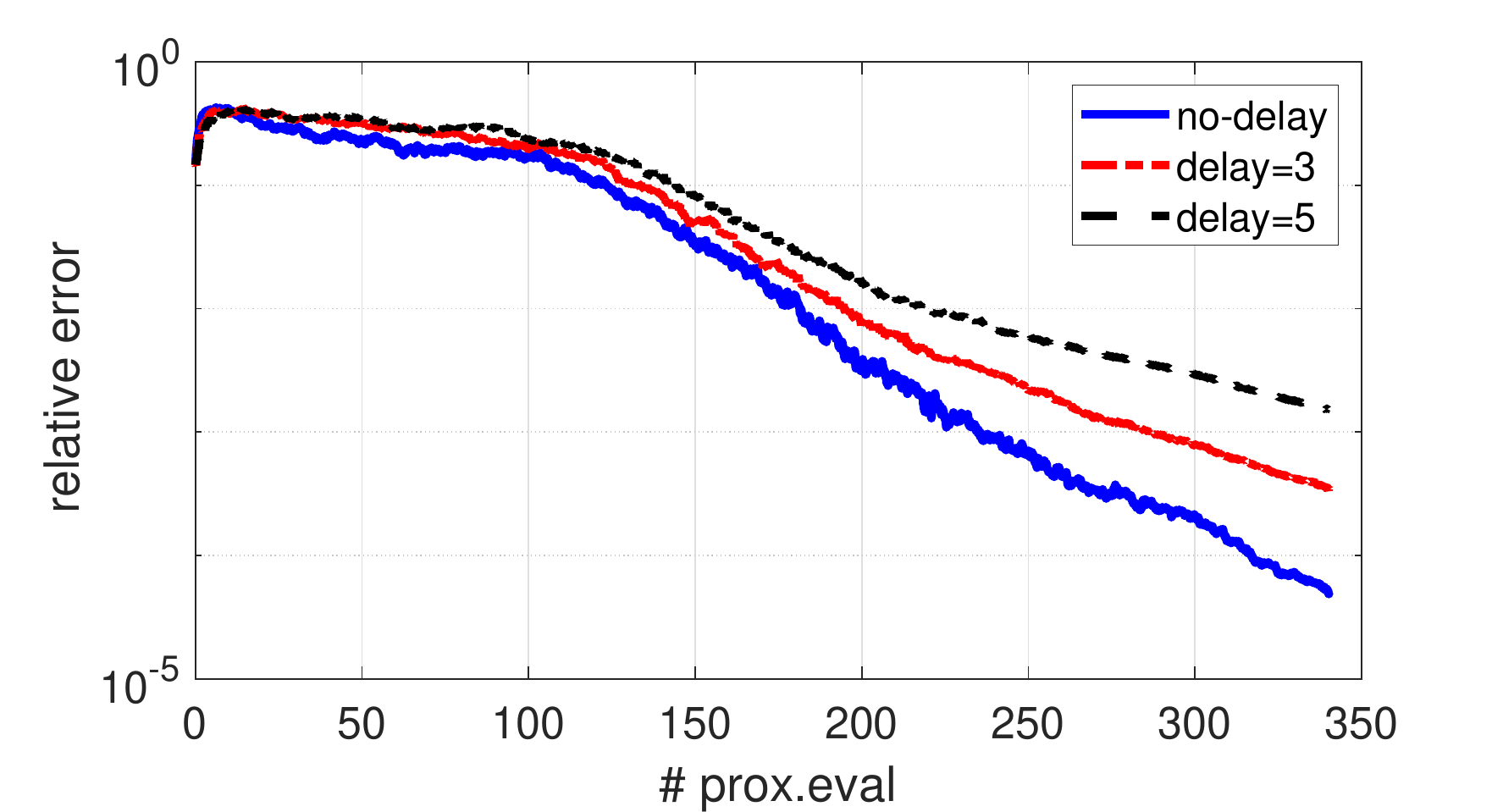}
    \end{minipage}%
    \begin{minipage}{0.46\textwidth}
        \centering
    \includegraphics[width=2.8in]{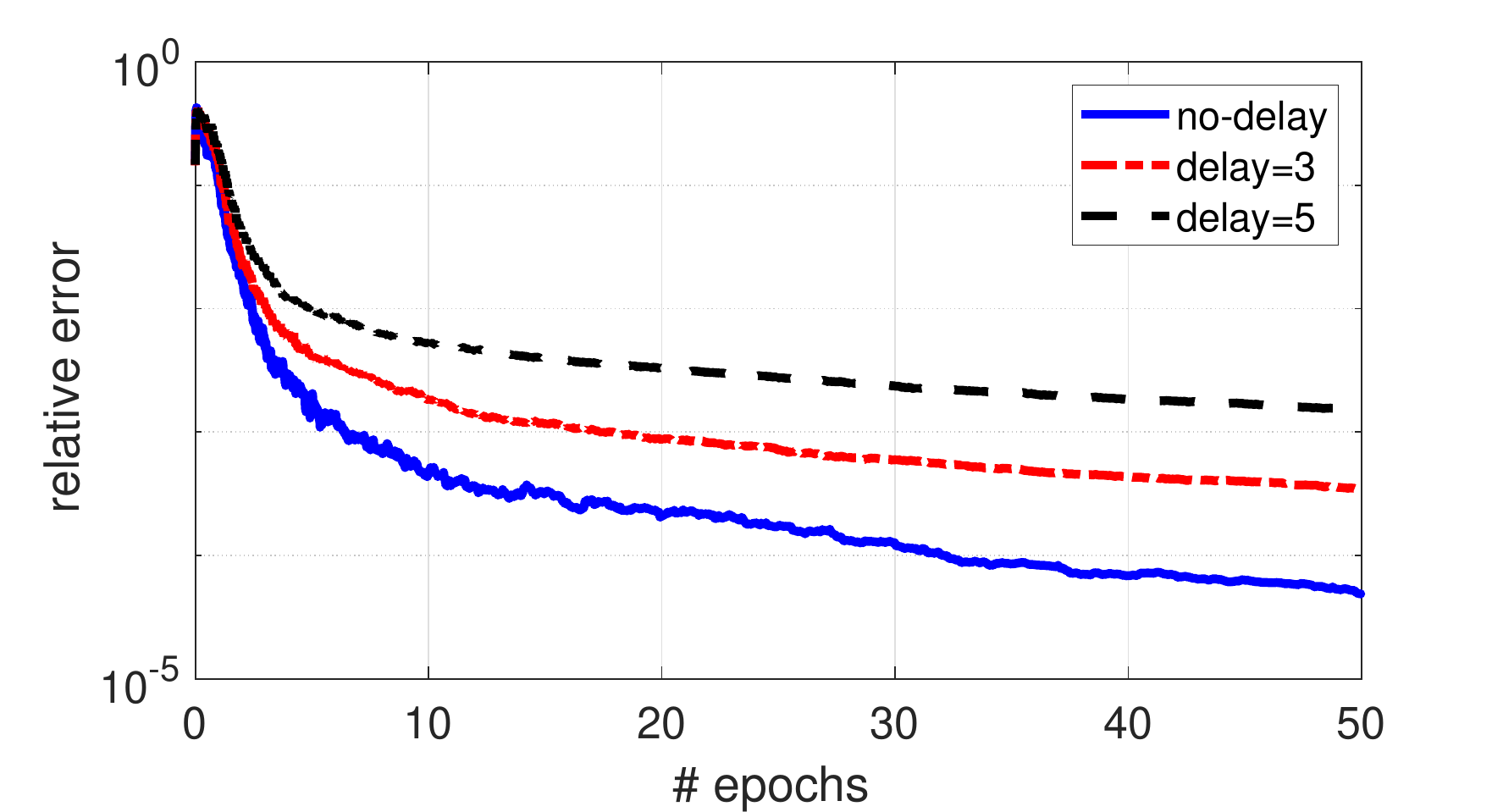}
    \end{minipage}
    \caption{ Performance of   Algorithm\ref{alg-prox-gradient} with Delays and without Delays   } \label{fig_delay}
\end{figure}

  \subsection{Nonlinear least squares}

We consider a binary classification problem on a data set
$\{x_i,y_i\}_{i=1}^N$, where $x_i\in \mathbb{R}^d$ and $y_i\in  \{0,1\}$ are
the $i$th feature vector and the corresponding label,  respectively.  We
consider  the  minimization of   empirical error: $ \min\limits_{w,b} {1\over
2N}\sum_{i=1}^N  \left(y_i-\phi(w^Tx_i+b) \right)^2, $ where $\phi(z)={1\over
1+e^{-z}}$ is the sigmoid function. We   apply   Algorithm
~\ref{alg-prox-gradient} to {\bf gisette}  from   LIBSVM library\footnotemark
\footnotetext{The data set is from
https://www.csie.ntu.edu.tw/~cjlin/libsvmtools/datasets/}  and investigate how
 batch-sizes    influence the  training loss and misclassification rate.
We partition the vector $w\in \mathbb{R}^d$ into $n=10$ blocks.   We implement
Algorithm \ref{alg-prox-gradient} with  $\alpha=0.2$,  where the batch-sizes
are set to be the constant  batch-sizes $N_i(k)\equiv 0.02N,  0.05N$   and  the
increasing batch-sizes $N_i(k)= \Gamma_i(k) ,$ $  (\Gamma_i(k))^2  $.  From Figure
\ref{fig-nls2},  we conclude that smaller  batch-sizes would  lead to  better
performance if we run the algorithm with  a relatively  small   amount of
samples (e.g, $N$); the mini-batch schemes may not perform well if the
batch-size is not suitably selected,  for instance, $0.05N$.  Favorable
behavior follows if the batch-size increases at a suitable rate, e.g.,
linearly.

\begin{figure}[!htb]
    \centering
    \begin{minipage}{.49\textwidth}
        \centering
      \includegraphics[width=2.8in]{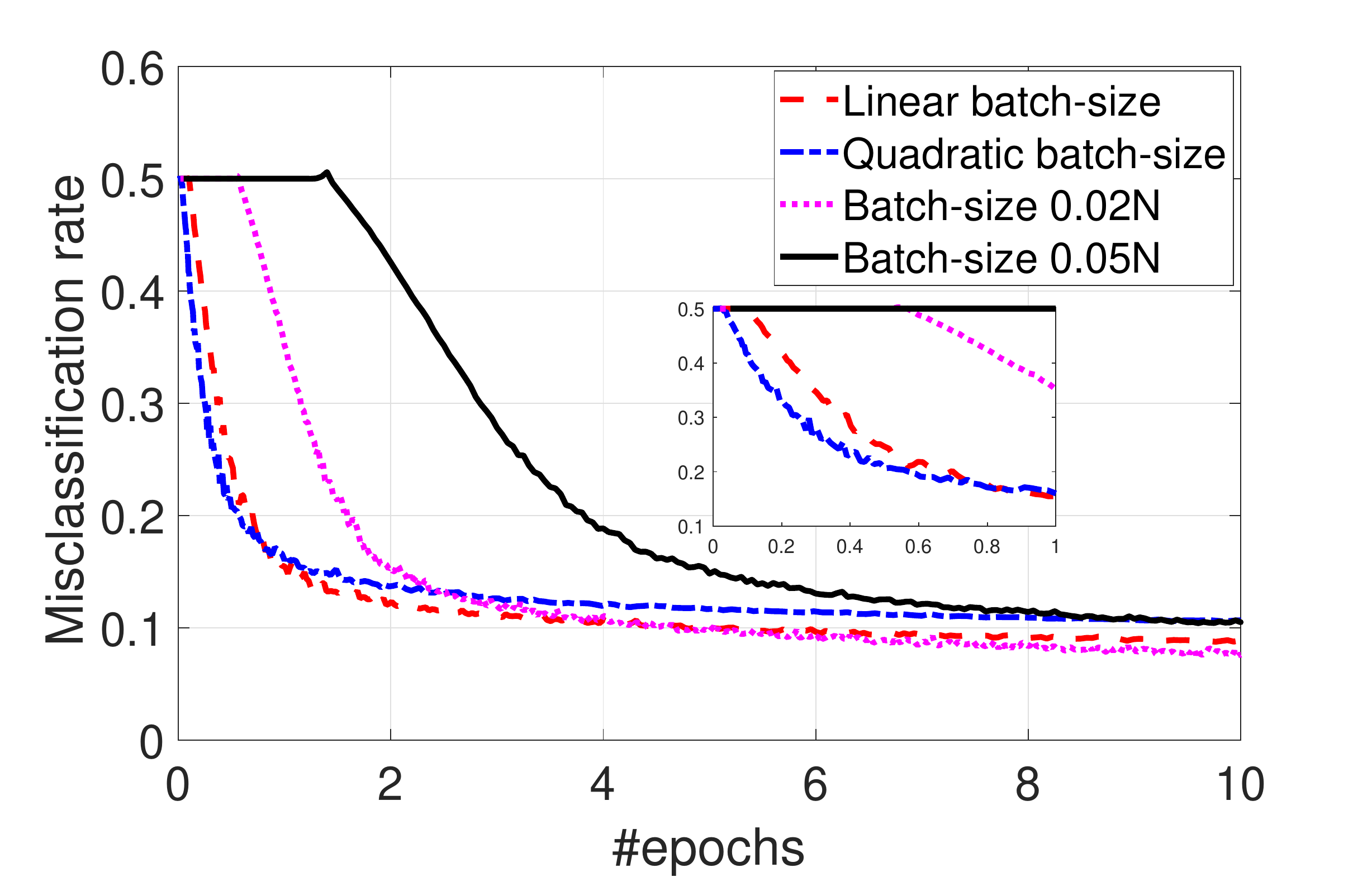}
    \end{minipage}%
    \begin{minipage}{0.49\textwidth}
        \centering
    \includegraphics[width=2.8in]{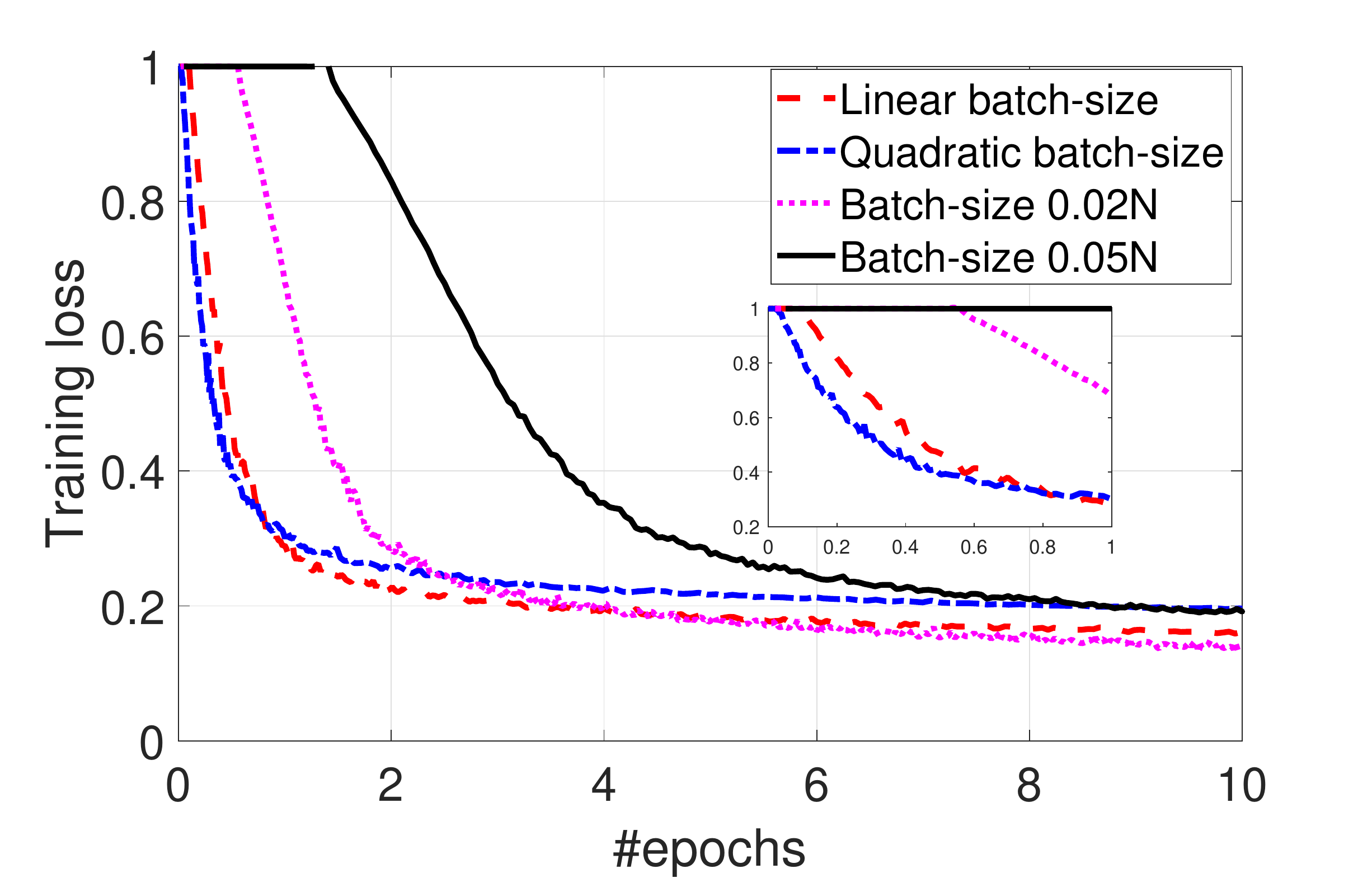}
    \end{minipage}
    \caption{Performance of Algorithm  \ref{alg-prox-gradient} with different batchsizes applied to the classification problem} \label{fig-nls2}
\end{figure} 

\section{Concluding remarks}\label{sec:conclusion}
 Existing block-based techniques for stochastic nonconvex optimization rely
on centrally mandated batch-sizes and steplengths bounded by the global
Lipschitz constant, leading to larger oracle complexities and poorer
performance (because of shorter steps), as well as higher informational  coordination requirements.
We consider  minimizing the sum of an expectation-valued smooth   nonconvex
function and a   nonsmooth separable   convex function through a {\em limited
coordination}  asynchronous variance-reduced method,  {reliant on
block-specific steplengths and random  decentralized batch-sizes.  The  almost sure convergence of the generated  iterates is established.  In addition, the scheme achieves the deterministic rate of $\mathcal{O}(1/K)$ with the rate and oracle complexities dependent on $L_{\rm ave}$ rather than $L_{\max}$. Furthermore, under the proximal PL requirement, the iterates provably  converge linearly (polynomially)  to the global optimum  in a mean sense when  batch-sizes grow
geometrically (polynomially). Notably, despite using randomized batch-size
sequences, we show that the deterministic iteration complexities may be
achieved.   
Specifically,  the schemes  achieve  the optimal oracle complexity when the considered  problem  admits a single block.
Finally, numerical studies  are  carried out to  support the theoretical findings  and reveal that schemes  leveraging  block-specific Lipschitz constants lead to significantly superior empirical behavior.

\appendix
\section{Proof of Lemma \ref{lem-sec3-2}}
Let $x_{-i}$ denote the coordinates of the variable $x$ except those correspond to block $i.$
 By   applying  Lemma \ref{pre-lem} to  the function  $\bar{f}\left( \cdot, x_{-i_k}(k) \right)$ and  Eqn. \eqref{def-xbar}  with $y=\bar{x}_{i_k}(k+1), z=x=x_{i_k  }(k), $ and $g= \nabla_{x_{i_k}} \bar{f}(x(k)) $,
 we obtain the following  inequality:
   \begin{align*}
& \bar{f}\left(x_{-i_k}(k),\bar{x}_{i_k}(k+1)\right)+r_{i_k}(\bar{x}_{i_k}(k+1)) \\&\leq  \bar{f}(x(k))+r_{i_k}(x_{i_k  }(k))+\left(\frac{L_{i_k}}{2}-{1\over  \alpha_i}\right)\| \bar{x}_{i_k}(k+1)-x_{i_k}(k)\|^2.
\end{align*}
Define $\tilde{x}(k+1)$ as follows:
\begin{align}\label{def-tildex}
 \tilde{x}_{i_k}(k+1) \triangleq \bar{x}_{i_k}(k+1) {\rm~~ and   ~~}\tilde{x}_j(k+1) \triangleq  x_j(k)~\forall j\neq i_k  .
\end{align} Then  $r_j(\tilde{x}_j(k+1))=r_j(x_j(k))~\forall j\neq i_k $,
and hence   we obtain the following bound:
 \begin{equation}\label{xbar-bd1}
 F\left(\tilde{x}(k+1)\right) \leq F\left(x(k)\right)  + \left(\frac{L_{i_k}}{2}-{1\over  \alpha_{i_k}}\right)\| \bar{x}_{i_k}(k+1)-x_{i_k}(k)\|^2.
\end{equation}
By applying Lemma \ref{pre-lem} to     the function  $\bar{f}\left( \cdot, x_{-i_k}(k) \right)$ and the update \eqref{proxg-rw} with  $y=x_{i_k}(k+1), z=\bar{x}_{i_k}(k+1), x=x_{i_k  }(k),$ and $ g= \nabla_{x_{i_k}} \bar{f}(x(k)) +w_{i_k}(k+1)$, one  obtains \begin{equation*}
\begin{split}
&\bar{f}\left(x_{-i_k}(k),x_{i_k}(k+1)\right)+r_{i_k}(x_{i_k}(k+1))  \leq \bar{f}\left(x_{-i_k}(k),\bar{x}_{i_k}(k+1)\right)+r_{i_k}(\bar{x}_{i_k}(k+1))\\&
-(x_{i_k}(k+1)-\bar{x}_{i_k}(k+1))^T w_{i_k}(k+1)  +\left(\frac{L_{i_k}}{2}-{1\over 2 \alpha_{i_k}}\right) \| x_{i_k  }(k+1)-x_{i_k  }(k)\|^2
\\& +\left(\frac{L_{i_k}}{2}+{1\over 2 \alpha_{i_k}}\right) \| \bar{x}_{i_k  }(k+1)-x_{i_k  }(k)\|^2- {1\over 2 \alpha_{i_k}} \| x_{i_k  }(k+1)-\bar{x}_{i_k}(k+1)\|^2.
\end{split}
\end{equation*}
Note that  $ \tilde{x}_j(k+1)=x_j(k)~\forall j\neq i_k$  and $ \tilde{x}_{i_k}(k+1) =\bar{x}_{i_k}(k+1)$ by definition \eqref{def-tildex}. Then by the definition of   $F(\cdot)$, we have the following
 \begin{equation} \label{x-bd}
 \begin{split} &F(x(k+1))  \leq   F(\tilde{x}(k+1))  -(x_{i_k}(k+1)-\bar{x}_{i_k}(k+1))^T w_{i_k}(k+1)
  \\&  +\left(\frac{L_{i_k}}{2}-{1\over 2 \alpha_{i_k}}\right) \| x_{i_k  }(k+1)-x_{i_k  }(k)\|^2 - {1\over 2 \alpha_{i_k}} \| x_{i_k  }(k+1)-\bar{x}_{i_k}(k+1)\|^2
\\&+\left(\frac{L_{i_k}}{2}+{1\over 2 \alpha_{i_k}}\right) \| \bar{x}_{i_k  }(k+1)-x_{i_k  }(k)\|^2.
 \end{split}
 \end{equation}  By  recalling that $-a^Tb\leq {1\over 2\alpha} \|a\|^2+{\alpha \over 2} \|b\|^2$, the following holds:
\begin{equation}\label{noise-bd}
\begin{split} & -(x_{i_k}(k+1)-\bar{x}_{i_k}(k+1))^T w_{i_k}(k+1) \\&  \leq {1\over 2 \alpha_{i_k}} \| x_{i_k}(k+1)-\bar{x}_{i_k}(k+1)\|^2+ {\alpha_{i_k} \over 2}   \|w_{i_k}(k+1) ||^2.
\end{split}
\end{equation}
Therefore,   by substituting  \eqref{noise-bd} into  \eqref{x-bd}, we obtain the following bound:
 \begin{equation}\label{x-bd0}
\begin{split}  F(x(k+1)) & \leq   F(\tilde{x}(k+1))    +\left(\frac{L_{i_k}}{2}-{1\over 2 \alpha_{i_k}}\right) \| x_{i_k}(k+1)-x_{i_k  }(k)\|^2
\\& +\left(\frac{L_{i_k}}{2}+{1\over 2 \alpha_{i_k}}\right)  \| \bar{x}_{i_k  }(k+1)-x_{i_k  }(k)\|^2  + {\alpha_{i_k} \over 2}   \|w_{i_k}(k+1)||^2 .
\end{split}
\end{equation}
   By adding inequalities  \eqref{xbar-bd1} and \eqref{x-bd0},
\begin{equation}\label{recursion-bd0}
\begin{split}
  F(x(k+1))  & \leq F\left(x(k)\right)  +\left(\frac{L_{i_k}}{2}-{1\over 2 \alpha_{i_k}}\right) \| x_{i_k}(k+1)-x_{i_k  }(k)\|^2
  \\& + {\alpha_{i_k} \over 2}   \|w_{i_k}(k+1)||^2+\left( L_{i_k }-{1\over 2 \alpha_{i_k}}\right)  \| \bar{x}_{i_k  }(k+1)-x_{i_k  }(k)\|^2.
\end{split}
\end{equation}
Note that for all $ i=1,\cdots,n,$  $\frac{L_i}{2}-{1\over 2 \alpha_i}\leq 0 $ by $\alpha_i \leq {1\over L_i}$. Then
  the second term on the right-hand side of  Eqn.  \eqref{recursion-bd0} is nonpositive, hence we can take out  this term from the upper bound of $  F(x(k+1)) $.  Since $x(k)$ is adapted to $\mathcal{F}_k$, by  taking expectations conditioned on $\mathcal{F}_k$ on both sides of \eqref{recursion-bd0},  we obtain that
 \begin{equation}\label{x-bd1}
\begin{split}
\mathbb{E}\big[ F(x(k+1)) |\mathcal{F}_k \big]&\leq  F(x(k))
 + \mathbb{E}\left[ \left(L_{i_k}-{1\over 2 \alpha_{i_k}}\right)\|\bar{x}_{i_k}(k+1)-x_{i_k  }(k)\|^2 |\mathcal{F}_k\right]
\\& +  {1 \over 2}     \mathbb{E}\big[   \alpha_{i_k}    \|w_{i_k}(k+1)||^2|\mathcal{F}_k\big].
\end{split}
\end{equation}
Note that  for  any $i \in \mathcal{N}$, $\bar{x}_i(k+1)$  is adapted to $\mathcal{F}_k$ by the definition \eqref{def-xbar}, and $i_k$ is independent of $\mathcal{F}_k$ by   Assumption \ref{ass-noise}(ii). Therefore,
by    \cite[Corollary 7.1.2 ]{chow2012probability}\footnotemark  \footnotetext{Let the random vectors
$X\in \mathbb{R}^m $ and $Y\in \mathbb{R}^n$ on $(\Omega, \mathcal{F},\mathbb{P})$
be independent of one another and let $f$ be a Borel function on $\mathbb{R}^{m \times n}$ with
$| \mathbb{E}[f(X,Y)] | \leq \infty$.  If for any $x\in  \mathbb{R}^m$,  $g(x)=\ \mathbb{E}[f(x,Y)]
 {\rm~ if~}| \mathbb{E}[f(x,Y)]|\leq \infty $ and  $g(x)=0~ {\rm otherwise}  $,  then $g$ is a Borel function with $g(X)=\mathbb{E}[f(X,Y)| \sigma(X)]$.} and $\mathbb{P}(i_k=i)=p_i$,   the following   holds a.s.:
	    \begin{align}\label{cond-expect}
 & \mathbb{E}\left[ \left(L_{i_k}-{1\over 2 \alpha_{i_k}}\right)\|\bar{x}_{i_k}(k+1)-x_{i_k  }(k)\|^2 |\mathcal{F}_k \right]
 =\sum_{i=1}^n p_i \left(L_i-{1\over 2 \alpha_i}\right)  \|\bar{x}_i(k+1)-x_i(k)\|^2.
\end{align}
Then   by substituting \eqref{cond-expect}   into \eqref{x-bd1},   we  obtain Eqn. \eqref{lem-x-bd}.
\hfill $\Box$

\section{Proof of Lemma \ref{sec4-lem1}.} By recalling that the gradient map $\nabla_{x_i} \bar{f}(x)$ is   $L_i $-Lipschitz continuous from Assumption \ref{ass-fun}(ii) and that  $\tilde{x}_j(k+1)=x_j(k)~\forall j\neq i_k $ by definition \eqref{def-tildex}, we have the following inequality:
\begin{equation*}
\bar{f}(\tilde{x}(k+1))\leq \bar{ f}(x(k))+(\tilde{x}_{i_k}(k+1)-x_{i_k  }(k))^T\nabla_{x_{i_k}}\bar{ f}(x(k))+{L_{i_k} \over 2}\| \tilde{x}_{i_k}(k+1)-x_{i_k  }(k)\|^2.
\end{equation*}
Using the definition of $\tilde{x}_{k+1}$ in \eqref{def-tildex},   we have that $r_j\left(\tilde{x}_j(k+1)\right)=r_j(x_j(k))~\forall j\neq i_k$,  $\tilde{x}_{i_k}(k+1)=\bar{x}_{i_k}(k+1)$. Thus, the the following relation holds
  \begin{align}\label{PL-bd1}
 F(\tilde{x}(k+1)) &\leq  F(x(k))+(\tilde{x}_{i_k}(k+1)-x_{i_k  }(k))^T\nabla_{x_{i_k}}\bar{ f}(x(k))+{L_{i_k}  \over 2}\| \tilde{x}_{i_k}(k+1)-x_{i_k  }(k)\|^2 \notag  \\& \quad +r_{i_k}(\tilde{x}_{i_k  }(k+1))-r_{i_k}(x_{i_k  }(k)) \notag
\\&\leq F(x(k))+(\bar{x}_{i_k}(k+1)-x_{i_k  }(k))^T\nabla_{x_{i_k}}\bar{ f}(x(k))\notag
\\& \quad+{1\over 2\alpha_{i_k} }\| \bar{x}_{i_k}(k+1)-x_{i_k  }(k)\|^2+r_{i_k}(\bar{x}_{i_k  }(k+1))-r_{i_k}(x_{i_k  }(k)),
\end{align}
where the last inequality holds by $\alpha_i <1/L_i  ~\forall i\in \mathcal{N}.$  Since   for any  $ i\in \mathcal{N},$  $\bar{x}_i(k+1) $  is adapted to $\mathcal{F}_k$ by its definition \eqref{def-xbar}, and  $i_k$ is independent of $\mathcal{F}_k$. Then,   by    \cite[Corollary 7.1.2 ]{chow2012probability} and $\mathbb{P}(i_k=i)=p_i$,  we have that
\begin{small}
   \begin{align*}
    &\mathbb{E} \Big [ (\bar{x}_{i_k}(k+1)-x_{i_k  }(k))^T\nabla_{x_{i_k}}\bar{ f}(x(k))+{\| \bar{x}_{i_k}(k+1)-x_{i_k  }(k)\|^2/( 2\alpha_{i_k})}  +r_{i_k}(\bar{x}_{i_k  }(k+1))-r_{i_k}(x_{i_k  }(k)) \big| \mathcal{F}_k \Big ]
\\& =\sum_{i=1}^n p_i \Big ( (\bar{x}_i(k+1)-x_i(k))^T\nabla_{x_{i }}\bar{ f}(x(k))  +{1\over 2\alpha_{i }}\| \bar{x}_i(k+1)-x_i(k)\|^2+r_{i }(\bar{x}_i(k+1))-r_{i }(x_i(k)) \Big)
\\& {{\scriptstyle \eqref{def-xbar} }\atop =}\sum_{i=1}^n p_i  \min_{y_i \in \mathbb{R}^{d_i}}\Big[\nabla_{x_{i }} \bar{f}(x(k))^T(y_i-x_i(k))+{1\over 2\alpha_{i } } \|y_i-x_i(k)\|^2+r_{i}(y_i)-r_{i}(x_i(k)) \Big]
\\& \leq p_{\min}\sum_{i=1}^n \min_{y_i \in \mathbb{R}^{d_i}}\Big[\nabla_{x_{i }} \bar{f}(x(k))^T(y_i-x_i(k))+{1\over 2\alpha_{i } } \|y_i-x_i(k)\|^2+r_{i}(y_i)-r_{i}(x_i(k)) \Big],
\end{align*}
\end{small}
  where the last  inequality follows by $ \min_{y_i \in \mathbb{R}^{d_i}} [\nabla_{x_{i }} \bar{f}(x(k))^T(y_i-x_i(k))+{1\over 2\alpha_i } \|y_i-x_i(k)\|^2+r_{i}(y_i)-r_{i}(x_i(k))   ]\leq 0$. Then by  $\alpha_i^{-1} \leq \alpha_{\min} ^{-1} $ and Assumption \ref{ass-PL},  the above equation  is further bounded by
\begin{equation}\label{bd-proximal}
\begin{split}
  &    p_{\min}\sum_{i=1}^n \min_{y_i \in \mathbb{R}^{d_i}}\Big[\nabla_{x_{i }} \bar{f}(x(k))^T(y_i-x_i(k))+{1 \over 2 \alpha_{\min}  } \|y_i-x_i(k)\|^2+r_{i}(y_i)-r_{i}(x_i(k)) \Big]
\\&=-   {  p_{\min} \alpha_{\min}  \over 2  }	 D_r(x(k),\alpha_{\min} ^{-1})
  \leq - {\alpha_{\min}  p_{\min}  \over 2}	 D_r(x(k),L_{\max})
 \\&  \leq - {\alpha_{\min} \mu   p_{\min}  } \big(F(x(k))-F^*\big)  ,
\end{split}
\end{equation}where the first  inequality follows from \cite[Lemma 1]{karimi2016linear} since  $D_r(x,\cdot) $ is nonnegative and  nondecreasing   in $(0,\infty)$   and $\alpha_{\min}^{-1} \geq L_{\max}  $.
Then by  taking unconditional  expectations on  both sides of \eqref{PL-bd1} and using \eqref{bd-proximal},  we obtain that
\begin{equation}\label{xbar-bd2}
\begin{split}
\mathbb{E}\big[ F(\tilde{x}(k+1))  \big]&\leq \mathbb{E}\big[ F(x(k)) \big]- \alpha_{\min} \mu   p_{\min}  \mathbb{E} \big[F(x(k))-F^*\big] .
\end{split}
\end{equation}
By taking unconditional   expectations on both sides of \eqref{xbar-bd1}  and using $\mathbb{P}(i_k=i)=p_i$, we obtain
 \begin{equation}\label{xbar-expect-bd1}
 \mathbb{E}\big[ F\left(\tilde{x}(k+1)\right)   \big] \leq \mathbb{E}\big[ F(x(k)) \big]  + \sum_{i=1}^n p_i \left({ L_i \over 2}-{1\over \alpha_i}\right)  \|\bar{x}_i(k+1)-x_i(k)\|^2.
\end{equation}
Adding  $(1-\beta) \times $ \eqref{xbar-bd2}  to $\beta \times $ \eqref{xbar-expect-bd1}  with $\beta \in (0.5,1)$,
we  obtain the following inequality:
 \begin{equation}\label{xbar-bd3}
\begin{split}
\mathbb{E}\big[ F(\tilde{x}(k+1))  \big]  & \leq \mathbb{E}\big[ F(x(k)) \big] +\beta\sum_{i=1}^n p_i \left({ L_i \over 2}-{1\over \alpha_i}\right)  \|\bar{x}_i(k+1)-x_i(k)\|^2
\\& -{\alpha_{\min}(1-\beta) \mu p_{\min} }   \mathbb{E} \big[F(x(k))-F^*\big] .
\end{split}
\end{equation}
Using  $\alpha_i <{1\over L_i},$ Assumption \ref{ass-noise}, and $\mathbb{P}(i_k=i)=p_i$, and
 by taking unconditional  expectations on both sides of \eqref{x-bd0},  the following holds:
 \begin{equation}\label{x-expect-bd0}
\begin{split}  \mathbb{E}\big[ F(x_{k+1}) \big]   &\leq   \mathbb{E}\big[ F(\tilde{x}(k+1))  \big]
 \\&\quad+ \sum_{i=1}^n p_i \left({ L_i \over 2}+{1\over 2 \alpha_i}\right)  \|\bar{x}_i(k+1)-x_i(k)\|^2
  +   { \sigma^2 \over 2} \sum_{i=1}^n\alpha_i p_i \mathbb{E}\left[  N_i(k)^{-1}\right]   .
\end{split}
\end{equation}
Therefore,  by adding  inequality  \eqref{x-expect-bd0} to  \eqref{xbar-bd3}  yields the following bound:
\begin{align}\label{x-bd2}
 \mathbb{E}\big[ F(x(k+1))\big] & \leq \mathbb{E}\big[ F(x(k))\big]  +
\sum_{i=1}^n p_i \left(\frac{L_i(1+\beta)}{2}-{2\beta-1\over 2 \alpha_i}\right) \|\bar{x}_i(k+1)-x_i(k)\|^2 \notag
\\& -{\alpha_{\min}(1-\beta) \mu p_{\min} }   \mathbb{E} \big[F(x(k))-F^*\big]  +   { \sigma^2 \over 2}\sum_{i=1}^n \alpha_i p_i \mathbb{E}\left[  N_i(k)^{-1}\right]   .
\end{align}
By recalling that  $0< \alpha_i\leq {2\beta-1\over L_i(1+\beta)}$, we get $\frac{L_i(1+\beta)}{2}-{2\beta-1\over 2 \alpha_i}\leq 0.$ Thus, by  subtracting $F^*$ from both sides of \eqref{x-bd2}, we  obtain \eqref{x-bd3}.
\hfill $\Box$

\end{document}